\documentclass[a4paper,final]{article}

\usepackage{amsmath,amssymb,amsthm,amsfonts,geometry,listings,graphicx,showkeys}

\newtheorem{assumption}{Assumption}
\newtheorem{theorem}{Theorem}
\newtheorem{proposition}{Proposition}
\newtheorem{lemma}{Lemma}

\providecommand{\R}{{\ensuremath{\mathbb{R}}}}

\renewcommand{\P}{{\ensuremath{\mathbb{P}}}}

\newcommand{\X}[1]{\ensuremath{X_{#1}}}
\newcommand{\XY}[1]{\ensuremath{\hat{X}_{#1}^{M}}}
\newcommand{\Y}[1]{\ensuremath{Y_{#1}^{M}}}

\newcommand{\Z}[1]{\ensuremath{Z_{#1}^{M}}}
\newcommand{\TZ}[1]{\ensuremath{\hat{Z}_{#1}^{M}}}
\newcommand{\TF}[1]{\ensuremath{\Phi_{#1}^M}}
\newcommand{\TB}[1]{\ensuremath{\Psi_{#1}^M}}
\newcommand{\C}[0]{\ensuremath{C(K)}}
\newcommand{\K}[0]{\ensuremath{C(K)}}
\newcommand{\fl}[2][M]{ \ensuremath{ \kern-1.0pt \left\lfloor \kern-0.75pt #2 \kern-0.75pt \right\rfloor_{_{\kern-2.0pt #1}} \kern-1.0pt } }
\newcommand{\cl}[2][M]{ \ensuremath{ \kern-1.0pt \left\lceil \kern-0.75pt #2 \kern-0.75pt \right\rceil_{_{\kern-2.0pt #1}} \kern-1.0pt } }
\newcommand{\m}[2][M]{ \ensuremath{ \kern-1.0pt \left< \kern-0.75pt #2 \kern-0.75pt \right>_{_{\kern-2.0pt #1}} \kern-1.0pt } }
\newcommand{\cll}[1]{ \ensuremath{ \kern-1.0pt \left\lceil \kern-0.75pt #1 \kern-0.75pt \right\rceil \kern-1.0pt } }
\newcommand{\I}[2]{ \ensuremath{ I_{ (#1), #2 } } }

\title{An exponential Wagner-Platen type scheme for SPDEs}
\author{Sebastian Becker, Arnulf Jentzen and Peter E. Kloeden}

\begin{document}

\maketitle


\begin{abstract}
The strong numerical approximation of semilinear 
stochastic partial differential equations (SPDEs) 
driven by infinite dimensional Wiener processes is investigated. 
There are   a number of results in the literature that show that Euler-type approximation methods 
converge strongly, under suitable assumptions, to the exact solutions of such SPDEs with strong order 
$ 1/2 $ or at least with strong order $ 1 / 2 - \varepsilon $
where $ \varepsilon > 0 $ is arbitrarily small. 
Recent results extend these results and show
that 
Milstein-type approximation methods 
converge, under suitable assumptions, 
to the exact solutions 
of such SPDEs
with strong order 
$ 1 - \varepsilon $.
It has also been shown 
that splitting-up approximation methods 
converge,
under suitable assumptions, with
strong order
$ 1 $ to the exact solutions 
of such SPDEs.
In this article   an  exponential
Wagner-Platen type numerical approximation method 
for such SPDEs is proposed and shown 
to  converge, under suitable assumptions, with strong order $ 3/2 - \varepsilon $ to the exact solutions of such SPDEs.
\end{abstract}

\section{Introduction}
\label{sec:intro}

We investigate   the strong numerical approximation of semilinear 
stochastic partial differential equations (SPDEs) 
driven by infinite dimensional Wiener processes. 
To illustrate the results, we concentrate in this introductory section on the following
simple example SPDE.
Let 
$ 
  H = L^2( (0,1); \R ) 
$
be the $ \R $-Hilbert space of equivalence classes
of Lebesgue square integrable functions from
$ ( 0, 1 ) $ to $ \R $,
let $ A \colon D( A ) \subset H \to H $ 
be the Laplacian with Dirichlet boundary conditions,
let $ ( \Omega, \mathcal{F}, \P, ( \mathcal{F}_t )_{ t \in [0,T] } ) $
be a stochastic basis (see, e.g., Appendix~E in Pr\'{e}v\^{o}t \& R\"{o}ckner~\cite{pr07}),
let $ r \in ( 1 ,\infty) $, $ \xi \in D(A) $,
let $ W \colon [0,T] \times \Omega \to H $
be a standard $ ( \mathcal{F}_t )_{ t \in [0,T] } $-Wiener process
with the covariance operator
$ Q := A^{ - 1 } \in L( H ) $
and let $ X \colon [0,T] \times \Omega \to H $
be an adapted stochastic process with continuous sample paths
satisfying the stochastic heat equation with linear multiplicative noise 
\begin{equation}
\label{eq:SPDE.intro}
  d X_t( x ) 
  =
  \left[ 
    \tfrac{ \partial^2 }{ \partial x^2 } X_t( x )
  \right] 
  dt
  +
  X_t(x) \, 
  dW_t(x)
  ,
\qquad 
  X_t( 0 ) = X_t( 1 ) = 0
  ,
\qquad 
  X_0(x) = \xi(x) 
\end{equation}
for 
$
  (t,x) \in [0,T] \times H 
$.
The convergence result in Theorem~\ref{thm:main} below can also be 
applied to a much more general class of SPDEs with more general covariance
operators $ Q $ 
(see Section~\ref{sec:sec2} for details) but for simplicity of presentation
we restrict ourselves to the SPDE~\eqref{eq:SPDE.intro} in this introductory section. 
Our goal is then to compute a strong numerical approximation of the
SPDE~\eqref{eq:SPDE.intro}.

There are a good number of results in the literature that show that 
{\it Euler-type approximation methods for SPDEs}
(such as the linear-implicit Euler method or the exponential Euler method;
see, e.g., Section~3.3.1 in Da Prato et al.~\cite{pjr10} for an overview on different Euler type approximations methods
for SPDEs)
converge to the solution process $ X $ of the SPDE~\eqref{eq:SPDE.intro} with strong order 
$ 1 / 2 $ or at least with strong order $ 1/2 - \varepsilon 
$ where $ \varepsilon \in ( 0, 1 / 2 ) $ is arbitrarily small
(see, e.g., \cite{cn13,h02,h03a,k11,lt13}). 
Further references on numerical approximations for SPDEs can also be found in the overview articles Gy{\"o}ngy~\cite{g02} and 
Jentzen \& Kloeden~\cite{jk09d}.

Recent results extend the above mentioned results for Euler type approximation methods 
and prove that {\it Milstein-type approximation methods for SPDEs}
converge with strong order $ 1 - \varepsilon $  or $1$ 
to the solution process $ X $
of the SPDE~\eqref{eq:SPDE.intro}
(see, e.g., \cite{bl12,bl09,jr10b,rk13,clp09,ms06,wg11}).
An overview on Milstein-type approximation methods for SPDEs
can also be found in Section~3.3.2 in Da Prato et al.~\cite{pjr10}.
Beside Milstein-type approximation methods,
it has also been established in the literature 
that {\it splitting-up approximation methods for SPDEs}
converge with
strong order
$ 1 $ to the solution process $ X $ of the SPDE~\eqref{eq:SPDE.intro} 
(see, e.g., Gy{\"o}ngy \& Krylov~\cite{gk03a}).
Further references for splitting-up methods can also be found in the overview article
Gy{\"o}ngy~\cite{g02}.
Beside Milstein type methods and splitting-up methods,
the choice of {\it suitable non-uniform time discretizations} is
another approach for obtaining higher order strong convergence rates 
for SPDEs; see \cite{mr07b,mr07a,mrw07,mrw08} for details on such methods and their optimality.
Moreover, in the case of SPDEs with additive noise,
higher order strong convergence rates for SPDEs
can also be obtained by using suitable {\it linear functionals
of the noise process in the numerical scheme};
see, e.g., \cite{ck06,bkh12,d10,j08a,jkw09,jk09b,kb13,lt10,lt10a,mtac12,wgt13}
and, e.g., (135)--(141) in Da Prato et al.~\cite{pjr10} for an overview.
For instance, in \cite{jk09b} it is proved in the case of additive noise,
that the accelerated exponential Euler method converges, under suitable assumptions, with strong order $ 1 - \varepsilon$ to the exact solution of
the SPDE under consideration.  Furthermore, higher order strong temporal converge rates of stochastic
Taylor schemes for spectral Galerkin discretizations
of SPDEs driven by one dimensional Brownian motions are established in Grecksch \& Kloeden \cite{gk96} and Kloeden \& Shott \cite{ks01}.

Here we  introduce an 
{\it exponential Wagner-Platen type numerical approximation method 
for SPDEs} (see \eqref{eq:scheme} below) and 
in Theorem~\ref{thm:main} below we prove 
that this method converges
with strong order $ 3/2 - \varepsilon $ to the solution process 
$ X $ of the SPDE~\eqref{eq:SPDE.intro}.
Further details can be found in Section~\ref{sec:scheme}
below.

In Section~\ref{sec:sec2} the abstract general setting used 
in this article is described.
Section~\ref{sec:scheme} introduces the above mentioned 
exponential Wagner-Platen method.
In addition, in Section~\ref{sec:scheme} we establish in
Proposition~\ref{prop:well_defined}
an a priori moment bound for the exponential Wagner-Platen method
and present  a convergence analysis theorem,
Theorem~\ref{thm:main}, of the exponential Wagner-Platen method.
Furthermore, Lemmas~\ref{lem:noise}--\ref{lem:int} in Section~\ref{sec:scheme} illustrate how the exponential Wagner-Platen 
method in Section~\ref{sec:scheme} can be simulated.
The proofs of Proposition~\ref{prop:well_defined},
Theorem~\ref{thm:main}, Lemma~\ref{lem:noise2} and Lemma~\ref{lem:int}
are postponed to Section~\ref{sec:proofs}.

\section{Setting}
\label{sec:sec2}

Throughout this article suppose 
that the following setting is fulfilled. 
Let $ T \in (0, \infty) $, 
let $ \left( \Omega, \mathcal{F}, \mathbb{P} \right) $ be a probability space 
with a normal filtration 
$ ( \mathcal{F}_t )_{ t \in [ 0, T] } $,
let 
$ \left( H, \left< \cdot , \cdot \right>_H, \left\| \cdot \right\|_H \right) $ 
and 
$ \left( U, \left< \cdot , \cdot \right>_U, \left\| \cdot \right\|_U \right) $ 
be separable $ \mathbb{R} $-Hilbert spaces,
let $ Q \in L(U) $ be a trace class operator
and let
$ 
  \left( U_0, \left< \cdot , \cdot \right>_{ U_0 }, \left\| \cdot \right\|_{ U_0 } \right) 
$ 
be the $ \mathbb{R} $-Hilbert space 
given by
$ U_0 = Q^{ 1 / 2 }( U ) $ 
and
$ 
  \left< v, w \right>_{ U_0 } = \big< Q^{ - 1 / 2 } v, Q^{ - 1 / 2 } w \big>_U 
$ 
for all 
$ v, w \in U_0 $.

\begin{assumption}[Linear operator A]
\label{ass:semigroup}
Let $ \mathcal{I} $ be a finite or countable set,
let $ \left( \lambda_i \right)_{i \in \mathcal{I} } $ $ \subset $ $(0, \infty) $ be a family of real numbers 
with $ \inf_{ i \in \mathcal{I} } \lambda_i \in (0,\infty) $,
let $ (e_i)_{ i \in \mathcal{I} } $ be an orthonormal basis of $H$ and 
let $ A \colon D(A) \subset H \rightarrow H $ be a linear operator 
with 
$ 
  D(A) = 
  \left\{ 
    w \in H 
    \colon
    \sum_{ i \in \mathcal{I} } 
    | \lambda_i |^2 
    | \langle e_i, w \rangle_H |^2 < \infty 
  \right\} 
$
and 
$
  Av
  =
  \sum_{ i \in \mathcal{I} }
  -\lambda_i
  \left<
    e_i, v
  \right>_H
  e_i
$
for all $v \in D(A)$.
\end{assumption}

In the following we denote
by
$
  ( 
    H_r, \left< \cdot, \cdot \right>_{ H_r },
    \left\| \cdot \right\|_{ H_r }
  ) 
$,
$ r \in \mathbb{R} $,
the $ \mathbb{R} $-Hilbert spaces
given by
$
  H_r = D( ( - A )^r )
$
and
$
  \left< v, w \right>_{ H_r }
  =
  \left< ( - A )^r v, ( - A )^r w \right>_H
$
for all 
$ v, w \in H_r $
and all
$ r \in \mathbb{R} $.

\begin{assumption}[Drift term $F$]
\label{ass:drift}
Let $ \gamma \in [1, \tfrac{3}{2} ) $, 
$ \alpha \in ( \gamma - 1 , \gamma] $ 
and let 
$ F \in C^2( H, H ) $ 
be a globally Lipschitz continuous mapping with 
$ 
  F( H_{\alpha} ) \subseteq H_{\alpha} 
$ 
and 
$ 
  \sup_{ v \in H_{\alpha} }
  \frac{
    \| 
      F(v) 
    \|_{ H_{\alpha} }
  }{
    1 + \| v \|_{ H_{\alpha} }
  }
  < 
  \infty
$.
\end{assumption}

\begin{assumption}[Diffusion term $B$]
\label{ass:diffusion}
Let $ \beta \in (\gamma -\tfrac{1}{2}, \gamma] $, $ \delta \in ( \gamma-1, \beta]  $ 
and let 
$ B \in C^2( H, HS(U_0,H) ) $
be a globally Lipschitz continuous mapping with $ B( H_{\beta} ) \subseteq HS( U_0,  H_{\beta} ) $,
$ 
  B'(v) \in 
  L( H_{\delta}, HS( U_0, H_{ \delta } ) 
  )
$ for all $ v \in H_{ \gamma } $
and
$
  \sup_{ v \in H_{\beta} } 
  \frac{
    \| 
      B(v) 
    \|_{ HS( U_0, H_{\beta} ) }
  }{
    1 + \| v \|_{ H_{\beta} }
  }
  +
  \sup_{ v \in H_{ \gamma } } 
  \left\| 
    B'(v) 
  \right\|_{ L( H_{\delta}, HS( U_0, H_{\delta} ) ) }
  <
  \infty
$.
\end{assumption}

\begin{assumption}[Initial value $\xi$]
\label{ass:initial}
Let $ \xi \colon \Omega \rightarrow H_{\gamma} $ be an $ \mathcal{F} / \mathcal{B}(H_{\gamma}) $-measurable mapping.
\end{assumption}

It is well known (see, e.g., 
Theorem~7.4 (i) in Da Prato \& Zabczyk~\cite{dz92})
that the above assumptions ensure the existence of an up to modifications unique predictable stochastic process 
$ 
  X \colon [0,T] \times \Omega \rightarrow H 
$ 
satisfying 
$ \int_0^T \big\| \X{s} \big\|_H^2 \, ds < \infty $ $ \mathbb{P} $-a.s.\ and
\begin{equation}
\label{eq:mild_solution}
  \X{t}
  =
  e^{ At } \xi
  +
  \int_0^t
  e^{ A(t-s) }
  F( \X{s} ) \, ds
  +
  \int_0^t
  e^{ A(t-s) }
  B( \X{s} ) \, dW_s
\end{equation}
$ \mathbb{P} $-a.s.\ for all $ t \in [0,T] $.

\section{An exponential Wagner-Platen type scheme for SPDEs}
\label{sec:scheme}

This section introduces and analyzes an exponential Wagner-Platen type approximation scheme
for the SPDE~\eqref{eq:mild_solution}.
To formulate this scheme, 
let $ \mathcal{J} $ be a finite or countable set 
and let $ ( g_j )_{ j \in \mathcal{J} } \subset U_0 $ be 
an arbitrary orthonormal basis of $ U_0 $.
Then let 
$ 
  \Y{m} \colon \Omega \rightarrow H_{ \gamma } 
$, 
$ 
  m \in \{ 0,1,\ldots,M\} 
$, 
$ M \in \mathbb{N}, 
$ be 
$ \mathcal{F}/\mathcal{B}(H_{ \gamma }) $-measurable mappings satisfying $ \Y{0} = \xi $ and
\begin{align}
\begin{split}
\label{eq:scheme}
&
  \Y{m+1}
 =
  e^{ A\frac{T}{2M} }
  \Bigg\{
    e^{ A\frac{T}{2M} }
    \Y{m}
    +
    \frac{T}{M} F( \Y{m} )
    +
    \frac{T^2}{2 M^2}
    F'( \Y{m} )
    \left[
      A \Y{m}
      +
      F( \Y{m} )
    \right]
\\&\quad+
    F'( \Y{m} )
    \left(
      \int_{ \frac{ m T }{ M } }^{ \frac{ (m + 1) T }{ M } }
      \int_{ \frac{ m T }{ M } }^s
      B( \Y{m} ) \, dW_u \, ds
    \right)
    +
    \sum_{ j \in \mathcal{J} }
    \frac{T^2}{4 M^2}
    F''( \Y{m} )
    \left(
      B( \Y{m} ) g_j, B( \Y{m} ) g_j
    \right)
\\&\quad+
    \int_{ \frac{ m T }{ M } }^{ \frac{ (m + 1) T }{ M } }
    B( \Y{m} ) \, dW_s
    +
    A
    \left[
      \int_{ \frac{ m T }{ M } }^{ \frac{ (m + 1) T }{ M } }
      \int_{ \frac{ m T }{ M } }^s
      B( \Y{m} ) \, dW_u \, ds
      -
      \frac{T}{2 M}
      \int_{ \frac{ m T }{ M } }^{ \frac{ (m + 1) T }{ M } }
      B( \Y{m} ) \, dW_s
   \right]
\\&\quad+
    \frac{T}{M}
    \int_{ \frac{ m T }{ M } }^{ \frac{ (m + 1) T }{ M } }
    B'( \Y{m} )
    \left(
      A \Y{m} + F( \Y{m} )
    \right) dW_s
    +
    \int_{ \frac{ m T }{ M } }^{ \frac{ (m + 1) T }{ M } }
    B'( \Y{m} )
    \left(
      \int_{ \frac{ m T }{ M } }^s
      B( \Y{m} ) \, dW_u
    \right) dW_s
\\&\quad-
    \int_{ \frac{ m T }{ M } }^{ \frac{ (m + 1) T }{ M } }
    \int_{ \frac{ m T }{ M } }^s
    B'( \Y{m} )
    \left(
      A \Y{m} + F( \Y{m} )
    \right) dW_u \, ds
\\&\quad+
    \frac{1}{2}
    \int_{ \frac{ m T }{ M } }^{ \frac{ (m + 1) T }{ M } }
    B''( \Y{m} )
    \left(
      \int_{ \frac{ m T }{ M } }^s
      B( \Y{m} ) \, dW_u,
      \int_{ \frac{ m T }{ M } }^s
      B( \Y{m} ) \, dW_u
    \right) dW_s
\\&\quad+
    \int_{ \frac{ m T }{ M } }^{ \frac{ (m + 1) T }{ M } }
    B'( \Y{m} )
    \left[
      \int_{ \frac{ m T }{ M } }^s
      B'( \Y{m} )
      \left(
        \int_{ \frac{ m T }{ M } }^u
        B( \Y{m} ) \, dW_v
      \right) dW_u
    \right] dW_s
  \Bigg\}
\end{split}
\end{align}
$ \mathbb{P} $-a.s.\ for all $ m \in \{ 0, 1, \ldots, M-1 \} $ and all $ M \in \mathbb{N} $. 
The setting in Section~\ref{sec:sec2} ensures that
the random variables 
$ Y^M_m \colon \Omega \to H_{ \gamma } $,
$ m \in \{ 0, 1, \dots, M \} $,
$ M \in \mathbb{N} $,
do indeed exist.
In addition, observe that the identity
\begin{equation}
\label{eq:conditional_expectation}
\begin{split}
&
  \mathbb{E}\!\left[ 
    F''( \Y{k} )
    \left(
      \int_{ \frac{ k T }{ M } }^t
      B( \Y{k} ) \, dW_u,
      \int_{ \frac{ k T }{ M } }^t
      B( \Y{k} ) \, dW_u
    \right)
    \, \Big| \,
    \mathcal{F}_{ \frac{ k T}{ M } }
  \right]
\\ & =
  \sum_{ j \in \mathcal{J} }
  \int_{ \frac{ k T }{ M } }^t
  F''( \Y{k} )
  \left(
    B( \Y{k} ) g_j, B( \Y{k} ) g_j
  \right) du
\end{split}
\end{equation}
$ \mathbb{P} $-a.s.\ for all 
$ 
  t \in [ \frac{ k T }{ M }, \frac{ (k+1) T }{ M } ] 
$, 
all 
$ 
  k \in \{ 0, 1, \ldots, M-1 \} 
$ and all 
$ 
  M \in \mathbb{N} 
$ 
illustrates that \eqref{eq:scheme}
does not depend on the special choice 
of the orthonormal basis 
$ ( g_j )_{ j \in \mathcal{J} } $
of $ U_0 $.
The following proposition 
establishes an a priori moment bound 
for the numerical approximations
$
  Y^M_m
$,
$ m \in \{ 0, 1, \dots, M \} $,
$ M \in \mathbb{N} $.

\begin{proposition}
\label{prop:well_defined}
There exists a universal non-decreasing function $ C \colon [0, \infty) \rightarrow [0, \infty) $ such that if the setting in Section~\ref{sec:sec2} is fulfilled, if $ p \in [2, \infty) $ and if 
\begin{align}
\notag
 &
  K
 :=
  \left\| A^{-1} \right\|_{ L(H) }
  +
  \frac{ 1 }{
    \left( 1 - \max( \gamma-\alpha, 2(\gamma-\beta), 2(\gamma-\delta-\tfrac{1}{2}), \gamma-\tfrac{1}{2} ) \right)
  }
  +
  T
  +
  p
  +
  \sup_{ v \in H_{\alpha} } 
  \frac{
    \left\| 
      F(v) 
    \right\|_{ H_{\alpha} }
  }{
    1 + \left\| v \right\|_{ H_{\alpha} }
  }
\\&+
\notag
  \sup_{ v \in H } 
  \left\| 
    F'(v) 
  \right\|_{ L( H ) }
  +
  \sup_{ v \in H_{\beta} } 
  \frac{
    \left\| 
      B(v) 
    \right\|_{ HS( U_0, H_{\beta} ) }
  }{
    1 + \left\| v \right\|_{ H_{\beta} }
  }
  +
  \sup_{ v \in H } 
  \left\| 
    B'(v) 
  \right\|_{ L( H, HS( U_0, H ) ) }
  +
  \sup_{ v \in H } 
  \left\| 
    B'( v ) 
  \right\|_{ L( H_{\delta}, HS( U_0, H_{\delta} ) ) }
\\&+
  \sup_{ v \in H } 
  \frac{
    \big( 
      \left\| 
        F''(v) 
      \right\|_{ L^{(2)}( H, H ) }
      +
      \left\| 
        B''(v) 
      \right\|_{ L^{(2)}( H, HS( U_0, H ) ) }
    \big)  
    \left\|
      B(v)
    \right\|_{ HS( U_0, H ) }^2
  }{
    1 + \left\| v \right\|_H
  }
  < 
  \infty ,
\end{align}
\begin{equation}
\label{eq:well_defined}
  \text{then}
\qquad
  \sup_{M \in \mathbb{N} }
  \sup_{m \in \{0,1,\ldots,M\}}
  \| \Y{m} \|_{ L^p( \Omega; H_{\gamma} ) } 
  \leq
  \K
  \left(
    1
    +
    \| \xi \|_{ L^p( \Omega; H_{\gamma} ) } 
  \right) 
  .
\end{equation}
\end{proposition}

The proof of Proposition~\ref{prop:well_defined} is postponed
to Subsection~\ref{sec:proof_of_prop1} below.
The next theorem estimates the strong temporal approximation error
of the numerical approximations
$
  Y^M_m
$,
$ m \in \{ 0, 1, \dots, M \} $,
$ M \in \mathbb{N} $.

\begin{theorem}
\label{thm:main}
There exists a universal non-decreasing function $ C \colon [0, \infty) \rightarrow [0, \infty) $ such that if the setting in Section~\ref{sec:sec2} is fulfilled and if
\begin{align}
\begin{split}
  &
  K
 :=
  \sup_{ m \in \{0,1,\ldots,M\}, M \in \mathbb{N} }
  \left\| \Y{m} \right\|_{ L^6( \Omega; H_{\gamma} ) } 
  +
  \left\| A^{-1} \right\|_{ L(H) }
  +
  T
  +
  \frac{ 1 }{ 
    \left( \min( \alpha+1, \beta+\tfrac{1}{2}, \delta+1, \tfrac{3}{2} ) - \gamma \right)
  }
\\&+
  \sup_{ v \in H_{\alpha} } 
  \frac{
    \left\| 
      F(v) 
    \right\|_{ H_{\alpha} }
  }{
    1 + \left\| v \right\|_{ H_{\alpha} }
  }
  +
  \sup_{ v \in H_{\beta} } 
  \frac{
    \left\| 
      B(v) 
    \right\|_{ HS( U_0, H_{\beta} ) }
  }{
    1 + \left\| v \right\|_{ H_{\beta} }
  }
  +
  \sup_{ v \in H_{ \gamma } } 
  \left\| 
    B'( v ) 
  \right\|_{ L( H_{\delta}, HS( U_0, H_{\delta} ) ) }
\\&+
  \sum_{ i=0 }^{ 2 }
  \sup_{ 
    v,w \in H , v \neq w 
  }
    \left[
  \tfrac{ 
      \left\|
        F^{(i)}(v)
        -
        F^{(i)}(w)
      \right\|_{ L^{(i)}(H,H) }
      +
      \left\|
        B^{(i)}(v)
        -
        B^{(i)}(w)
      \right\|_{ L^{(i)}(H,HS( U_0, H )) }
  }{ 
    \left\| v - w \right\|_H
  }
    \right]
  <
  \infty ,
\end{split}
\end{align}
\begin{equation}
  \text{then it holds for all $ M \in \mathbb{N} $ that}
  \qquad
  \sup_{ m \in \{ 0,1,\ldots,M\} } 
  \left\|
    X_{mh}
    -
    \Y{m}
  \right\|_{ L^2( \Omega; H ) }
  \leq
  \C 
  \, 
  M^{-\gamma }
  .
\end{equation}
\end{theorem}

The proof of Theorem~\ref{thm:main} is given
in Subsection~\ref{sec:proof_of_thm1} below.
The following lemmas (Lemma~\ref{lem:noise}, 
Lemma~\ref{lem:noise2} and Lemma~\ref{lem:int})
show under suitable assumptions how 
the scheme \eqref{eq:scheme} can be simulated.
Lemma~\ref{lem:noise} is a slightly more general
statement than display (83) in \cite{jr10b} 
(see also Remark~1 and Subsection~5.7 in \cite{jr10b}).

\begin{lemma}[Commutative noise of the first kind for SPDEs]
\label{lem:noise}
Assume the setting in Section~\ref{sec:sec2} and assume for all $ v \in H $ that the bilinear Hilbert Schmidt operator $ B'(v)(B(v)) \in H^{(2)}( U_0, H ) $ is symmetric, i.e., assume that
\begin{equation}
\label{eq:comm}
  \Big( 
    B'(v)
    \big(
      B(v) u_1
    \big)
  \Big) (u_2)
  =
  \Big( 
    B'(v)
    \big(
      B(v) u_2
    \big)
  \Big) (u_1)
\end{equation}
for all $ v\in H $ and all $ u_1, u_2 \in U_0 $. If $ U = U_0 $ (which is equivalent to $ dim(U) < \infty$), then
\begin{align}
\label{eq:b'b} 
\begin{split}
 &\int_{ t_0 }^{ t }
  B'( Z )
  \left( 
    \int_{ t_0 }^s
    B( Z ) \, dW_u 
  \right) dW_s
\\&=
  \frac{1}{2}
  B'( Z )
  \Big( 
    B( Z )
    \left(W_t - W_{ t_0 }\right)
  \Big)
  \left(W_t - W_{ t_0 }\right)
  -
  \frac{\left( t - t_0 \right)}{2}
  \sum_{ i \in \mathcal{J} }
  B'( Z )
  \Big( 
    B( Z ) \, g_i
  \Big) \, g_i
\end{split}
\end{align}
$ \mathbb{P} $-a.s. for all $ \mathcal{F}_{t_0} / \mathcal{B}(H) $-measurable mappings $ Z \colon \Omega \rightarrow H $ and all $ t_0, t \in [0,T] $ with $ t_0 \leq t $.
\end{lemma}

The proof of Lemma~\ref{lem:noise} is entirely analogous 
to the proof of (83) in \cite{jr10b}
and therefore omitted.
The next lemma treats the case of \emph{commutative noise of 
second kind for SPDEs}.
Assumption~\eqref{eq:comm2}
is the abstract coordinate free analogue 
of (4.13) in Section~10.4 in Kloeden \&\ Platen~\cite{kp92}.



\begin{lemma}[Commutative noise of the second kind for SPDEs]
\label{lem:noise2}
Assume the setting in Section~\ref{sec:sec2} and assume for all $ v \in H $ that the trilinear Hilbert Schmidt operators $ B'(v)(B'(v)(B(v))) \in H^{(3)}( U_0, H ) $ and $ B''(v)(B(v),B(v)) \in H^{(3)}( U_0, H ) $ are symmetric, i.e., assume that
\begin{align}
\label{eq:comm2}
\begin{split}
 &\Big( 
    B'(v)
    \Big(
      B'(v)
      \big(
        B(v) u_1
      \big)
    \Big) u_2
    +
    B''(v)
    \big( B(v) u_1, B(v) u_2 \big)
  \Big) (u_3)
\\&=
  \Big( 
    B'(v)
    \Big(
      B'(v)
      \big(
        B(v) u_{\pi(1)}
      \big)
    \Big) u_{\pi(2)}
    +
    B''(v)
    \big( B(v) u_{\pi(1)}, B(v) u_{\pi(2)} \big)
  \Big) (u_{\pi(3)})
\end{split}
\end{align}
for all $ v\in H $, all $ u_1, u_2, u_3 \in U_0 $ and all $ \pi \in S_3 $. If $ U = U_0 $ (which is equivalent to $ dim(U) < \infty$), then
\begin{align}
\notag
  &\int_{ t_0 }^{ t }
  B'( Z )
  \left( 
    \int_{ t_0 }^s
    B'( Z ) 
    \left( 
      \int_{ t_0 }^u
      B( Z ) \, dW_v 
    \right) dW_u 
  \right) dW_s
  +
  \frac{1}{2}
  \int_{ t_0 }^{ t }
  B''( Z )
  \left(
    \int_{ t_0 }^s
    B( Z ) \, dW_u,
    \int_{ t_0 }^s
    B( Z ) \, dW_u
  \right) dW_s
\\&=
\notag
  \frac{1}{6}
  B'( Z )
  \Big( 
    B'( Z ) 
    \Big( 
      B( Z ) \left(W_t - W_{t_0}\right)
    \Big) \left(W_t - W_{t_0}\right)
  \Big)
  \left(W_t - W_{t_0}\right)
\\&+
\label{eq:b'b'b}
  \frac{1}{6}
  B''( Z )
  \Big(
    B( Z ) \left(W_t - W_{t_0}\right),
    B( Z ) \left(W_t - W_{t_0}\right)
  \Big) 
  \left(W_t - W_{t_0}\right)
\\&-
\notag
  \frac{\left( t - t_0 \right)}{2}
  \sum_{ i \in \mathcal{J} }
  \left[ 
  B'( Z )
  \Big( 
    B'( Z ) 
    \Big( 
      B( Z ) \, g_i 
    \Big) g_i 
  \Big)
  \left(W_t - W_{t_0}\right)
  +
  B''( Z )
  \Big(
    B( Z ) \, g_i,
    B( Z ) \, g_i
  \Big) 
  \left(W_t - W_{t_0}\right)
  \right]
\end{align}
$ \mathbb{P} $-a.s. for all $ \mathcal{F}_{t_0} / \mathcal{B}(H) $-measurable mappings $ Z \colon \Omega \rightarrow H $ and all $ t_0, t \in [0,T] $ with $ t_0 \leq t $.
\end{lemma}

The proof of Lemma~\ref{lem:noise2}
is postponed to Subsection~\ref{sec:proof_of_lemnoise2} below.
Combining Lemma~\ref{lem:noise} and Lemma~\ref{lem:noise2} shows that if $ U = U_0 $ and if \eqref{eq:comm} and \eqref{eq:comm2} are fulfilled then the numerical scheme \eqref{eq:scheme} satisfies 
\begin{align}
\notag
&
  \Y{m+1}
 =
  e^{ A\frac{T}{2M} }
  \Bigg\{
    e^{ A\frac{T}{2M} }
    \Y{m}
    +
    \frac{T}{M} F( \Y{m} )
    +
    F'( \Y{m} )
    \frac{T^2}{2 M^2}
    \left[
      A \Y{m}
      +
      F( \Y{m} )
    \right]
\\&\quad+
\notag
    F'( \Y{m} )
    \left( 
    B( \Y{m} ) 
    \int_{\frac{mT}{M}}^{\frac{(m+1)T}{M}}
    \left( W_s - W_{\frac{mT}{M}} \right) ds
    \right)
    +
    \sum_{ j \in \mathcal{J} }
    \frac{T^2}{4 M^2}
    F''( \Y{m} )
    \left(
      B( \Y{m} ) g_j, B( \Y{m} ) g_j
    \right)
\\&\quad+
\notag
    B( \Y{m} )
    \left(W_{\frac{(m+1)T}{M}} - W_{\frac{mT}{M}}\right)
    +
    A
    \left[
      B( \Y{m} )
      \left( 
        \int_{\frac{mT}{M}}^{\frac{(m+1)T}{M}}
        \left( W_s - W_{\frac{mT}{M}} \right) ds
        -
        \frac{T}{2 M}
        \left(W_{\frac{(m+1)T}{M}} - W_{\frac{mT}{M}}\right)
      \right)
   \right]
\\&\quad+
\notag
    \Big(
      B'( \Y{m} )
      \left(
        A \Y{m} + F( \Y{m} )
      \right)
    \Big) \!
    \left(
      \frac{T}{M}
      \left(W_{\frac{(m+1)T}{M}} - W_{\frac{mT}{M}}\right)
      -
      \int_{\frac{mT}{M}}^{\frac{(m+1)T}{M}}
      \left( W_s - W_{\frac{mT}{M}} \right) ds
    \right)  
\\&\quad+
\label{eq:scheme_finite_noise}
    \frac{1}{2}
    B'( \Y{m} )
    \left( 
      B( \Y{m} ) 
      \left(W_{\frac{(m+1)T}{M}} - W_{\frac{mT}{M}}\right)
    \right) 
    \left(W_{\frac{(m+1)T}{M}} - W_{\frac{mT}{M}}\right)
\\&\quad+
\notag
    \frac{1}{6}
    B''( \Y{m} )
    \left(
      B( \Y{m} ) \left(W_{\frac{(m+1)T}{M}} - W_{\frac{mT}{M}}\right),
      B( \Y{m} ) \left(W_{\frac{(m+1)T}{M}} - W_{\frac{mT}{M}}\right)
    \right) 
    \left(W_{\frac{(m+1)T}{M}} - W_{\frac{mT}{M}}\right)
\\&\quad+
\notag
    \frac{1}{6}
    B'( \Y{m} )
    \left( 
      B'( \Y{m} )
      \left( 
        B( \Y{m} ) 
        \left(W_{\frac{(m+1)T}{M}} - W_{\frac{mT}{M}}\right)
      \right) \left(W_{\frac{(m+1)T}{M}} - W_{\frac{mT}{M}}\right)
    \right) \left(W_{\frac{(m+1)T}{M}} - W_{\frac{mT}{M}}\right)
\\&\quad-
\notag
    \frac{T}{2M}
    \sum_{ j \in \mathcal{J} }
    B'( \Y{m} )
    \left( 
      B( \Y{m} ) \, g_j
    \right) g_j
    -
    \frac{T}{2M}
    \sum_{ j \in \mathcal{J} }
    B''( \Y{m} )
    \left(
      B( \Y{m} ) \, g_j,
      B( \Y{m} ) \, g_j
    \right)  
    \left(W_{\frac{(m+1)T}{M}} - W_{\frac{mT}{M}}\right)
\\&\quad-
\notag
    \frac{T}{2M}
    \sum_{ j \in \mathcal{J} }
    B'( \Y{m} )
    \left( 
      B'( \Y{m} )
      \left(
        B( \Y{m} ) \, g_j
      \right) g_j
    \right) \left(W_{\frac{(m+1)T}{M}} - W_{\frac{mT}{M}}\right)
  \Bigg\}
\end{align}
$ \mathbb{P} $-a.s. for all $ m \in \{ 0,1, \ldots, M-1 \} $ and all $ M \in \mathbb{N} $.
The next lemma, Lemma~\ref{lem:int}, illustrates 
for all $ t_0, t \in [0,T] $ with $ t_0 \leq t $
how the Gaussian distributed random variable
$
  \big( W_t - W_{ t_0 } , \int_{ t_0 }^t ( W_s - W_{ t_0 } ) \, ds \big)
  \in U \times U
$
can be simulated. 
Lemma~\ref{lem:int} generalizes (4.2)--(4.3) in Section~10.4 in 
Kloeden \&\ Platen~\cite{kp92} for finite dimensional SODEs to infinite dimensional 
Wiener processes. The proof of Lemma~\ref{lem:int}
is given in Subsection~\ref{sec:proof_of_lemint} below.

\begin{lemma}[Covariance operator]
\label{lem:int}
Assume the setting in Section~\ref{sec:sec2}. Then
\begin{align}
\begin{split}
  \left(
  \operatorname{Cov}\!\left(
    \begin{array}{c}
        W_t - W_{t_0} \\
        \int_{t_0}^t \left( W_s - W_{t_0} \right) ds
    \end{array} 
  \right)
  \right)
  \left(
    \begin{array}{c}
      u_1 \\
      u_2
    \end{array}
  \right)
  =
  \left(
    \begin{array}{c}
      \left( t - t_0 \right) Q u_1 + \frac{1}{2} \left( t - t_0 \right)^2 Q u_2 \\
      \frac{1}{2} \left( t - t_0 \right)^2 Q u_1 + \frac{1}{3} \left( t - t_0 \right)^3 Q u_2 
    \end{array}
  \right)
\end{split}
\end{align}
for all $ u_1, u_2 \in U $ and all $ t_0, t \in [0,T] $ with $ t_0 \leq t $.
\end{lemma}

\section{Proofs}
\label{sec:proofs}
Throughout the rest of this article we use the mappings 
$ \fl{\cdot}, \cl{\cdot}, \m{\cdot} \colon [0,T] \rightarrow [0,T] $, $ M \in \mathbb{N} $, given by 
\begin{equation}
  \fl{t}
  :=
  \max\!\left\{
    s \in
    \left\{
      0, \tfrac{T}{M}, \ldots ,\tfrac{(M-1)T}{M}, T
    \right\}
    \colon
    s \leq t
  \right\}, \quad
  \cl{t}
  :=
  \min\!\left\{
    s \in
    \left\{
      0, \tfrac{T}{M}, \ldots ,\tfrac{(M-1)T}{M}, T 
    \right\}
    \colon
    s \geq t
  \right\}
\end{equation}
and $ \m{t} := \tfrac{1}{2} \left( \fl{t} + \cl{t} \right) $ for all $ t \in [0,T] $ and all $ M \in \mathbb{N} $. Moreover, let $ \XY{} \colon [0, T] \times \Omega \rightarrow H $, $ M \in \mathbb{N} $, be optional measurable stochastic processes satisfying 
$ 
  \sup_{ s \in [0,T] } 
  \mathbb{E}\big[ 
    \| \XY{s} \|_H^6
  \big] < \infty 
$ 
for all $ M \in \mathbb{N} $ and
\begin{align}
\label{eq:helper_process_1}
\begin{split}
  \XY{t}
 &=
  e^{ A(t-\fl{t}) } \Y{ \left( M / T \right) \fl{t} }
  +
  \int_{\fl{t}}^t
  e^{A(t-s)} F( \XY{s} ) \, ds
  +
  \int_{\fl{t}}^t
  e^{A(t-s)} B( \XY{s} ) \, dW_s
\end{split}
\end{align}
$ \mathbb{P} $-a.s.\ for all $ t \in [0, T] $ and all $ M \in \mathbb{N} $ (see, e.g., Theorem~7.4 (i) in Da Prato \& Zabczyk~\cite{dz92}). In addition, let $ \TF{} \colon [0, T] \times \Omega \rightarrow H $, $ M \in \mathbb{N} $, and $ \TB{} \colon [0, T] \times \Omega \rightarrow HS( U_0, H ) $, $ M \in \mathbb{N} $, be optional measurable stochastic processes satisfying
\begin{align}
\begin{split}
  \TF{t}
 &=
  F( \XY{\fl{t}} )
  +
  F'(\XY{\fl{t}})
  \int_{\fl{t}}^t
  \left(
    A \XY{\fl{t}}
    +
    F(\XY{\fl{t}}) 
  \right) du
  +
  F'(\XY{\fl{t}})
  \int_{\fl{t}}^t
  B(\XY{\fl{t}}) \, dW_u
\\&+
\label{eq:taylor_F}
  \frac{1}{2}
  \sum_{ j \in \mathcal{J} }
  \int_{\fl{t}}^t
  F''( \XY{\fl{t}} )
  \left(
    B( \XY{\fl{t}} ) g_j, B( \Y{\fl{t}} ) g_j
  \right) du
\end{split}
\end{align}
$ \mathbb{P} $-a.s.\ and
\begin{align}
\begin{split}
\label{eq:taylor_B}
  \TB{t}
 &=
  B(\XY{\fl{t}})
  +
  B'(\XY{\fl{t}})
  \int_{\fl{t}}^t
  \left(
    B(\XY{\fl{t}})
    +
    B'(\XY{\fl{t}})
    \int_{\fl{t}}^u
    B(\XY{\fl{t}}) \, dW_v
  \right) dW_u
\\&+
  B'(\XY{\fl{t}})
  \int_{\fl{t}}^t \!\!\!
  \left(
    A \XY{\fl{t}}
    +
    F(\XY{\fl{t}}) 
  \right) du
  +
  \frac{1}{2}
  B''(\XY{\fl{t}})
  \left(
    \int_{\fl{t}}^t \!\!\!
    B(\XY{\fl{t}}) \, dW_u,
    \int_{\fl{t}}^t \!\!\!
    B(\XY{\fl{t}}) \, dW_u
  \right)
\end{split}
\end{align}
$ \mathbb{P} $-a.s.\ for all $ t \in [0, T] $ and all $ M \in \mathbb{N} $. Note that It{\^o}'s formula shows
\begin{align}
\label{eq:ito_1}
\begin{split}
 &A
  \left[
    \int_{\frac{mT}{M}}^{ \frac{(m+1)T}{M} }
    \int_{\frac{mT}{M}}^s
    B(\Y{m}) \, dW_u \, ds
    -
    \frac{T}{2M}
    \int_{\frac{mT}{M}}^{ \frac{(m+1)T}{M} }
    B(\Y{m}) \, dW_s
  \right]
\\&=
  A
  \left[
    \frac{T}{2M}
    \int_{\frac{mT}{M}}^{ \frac{(m+1)T}{M} }
    B(\XY{\fl{t}}) \, dW_s
    -
    \int_{\frac{mT}{M}}^{ \frac{(m+1)T}{M} }
    \left(s-\tfrac{mT}{M} \right) 
    B(\Y{m}) \, dW_s
  \right]
\\&=
  \int_{\frac{mT}{M}}^{ \frac{(m+1)T}{M} }
  \left(
    \m{s}
    -
    s
  \right) \! A
  B( \XY{\fl{s}} ) \, dW_s
\end{split}
\end{align}
$ \mathbb{P} $-a.s.\ and
\begin{align}
\label{eq:ito_2}
\begin{split}
 &\frac{T}{M}
  \int_{\frac{mT}{M}}^{ \frac{(m+1)T}{M} }
  B'(\Y{m})
  \left(
    A \Y{m} + F( \Y{m} )
  \right) dW_s
  -
  \int_{\frac{mT}{M}}^{ \frac{(m+1)T}{M} }
  \int_{\frac{mT}{M}}^s
  B'(\Y{m})
  \left(
    A \Y{m} + F( \Y{m} )
  \right) dW_u \, ds
\\&=
  \int_{\frac{mT}{M}}^{ \frac{(m+1)T}{M} }
  \int_{\fl{s}}^s
  B'(\XY{\fl{s}})
  \left(
    A \XY{\fl{s}} + F( \XY{\fl{s}} )
  \right) du \, dW_s
\end{split}
\end{align}
$ \mathbb{P} $-a.s.\ for 
all 
$ m \in \{ 0, 1, \ldots, M-1 \} $ 
and all $ M \in \mathbb{N} $. 
Next we combine 
\eqref{eq:scheme}, \eqref{eq:ito_1} 
and \eqref{eq:ito_2} 
to obtain that
\begin{align}
\begin{split}
  \Y{m}
 &=
\label{eq:scheme_new}
  e^{ A\frac{mT}{M} }
  \xi
  +
  \int_{0}^{\frac{mT}{M}}
  e^{ A(\frac{mT}{M}-\m{s}) } \,
  \TF{s} \, ds
  +
  \int_{0}^{\frac{mT}{M}}
  e^{ A(\frac{mT}{M}-\m{s}) }
  \left( 
    \TB{s}
    +
    \left( \m{s} \!-\! s \right) \! A
    B( \XY{\fl{s}} )
  \right) dW_s
\end{split}
\end{align}
$ \mathbb{P} $-a.s.\ for all $ m \in \{ 0,1, \ldots, M \} $ and all $ M \in \mathbb{N} $.

\subsection{Proof of Proposition \ref{prop:well_defined}}
\label{sec:proof_of_prop1}

Throughout this proof 
$ C \colon [0,\infty) \to [0,\infty) $
is a universal non-decreasing
function which changes from line to line.
Let
$ \theta \in [0,\infty) $
be defined by
$ 
  \theta 
  := 
  \max( 
    \gamma - \alpha, 
    \gamma - \frac{ 1 }{ 2 }, 
    2 ( \gamma - \beta ), 
    2 ( \gamma - \delta - \frac{ 1 }{ 2 } ) 
  ) 
$. 
Observe that Assumption~\ref{ass:drift}
and Assumption~\ref{ass:diffusion}
ensure that $ \theta < 1 $.
Next note that
\begin{align}
\begin{split}
 &\left\| 
    e^{ \frac{A}{2}(t-\m{s}) }
    \TF{s}
  \right\|_{ L^p( \Omega; H_{(\gamma-\theta)} ) }^2
\\&\leq
  \K
  \left\|
    \left( -A \right)^{ (\gamma-\theta-\alpha) }
  \right\|_{ L(H) }^2
  \left\|
    F( \XY{\fl{s}} )
  \right\|_{ L^p( \Omega; H_{\alpha} ) }^2
\\&+
  \K
  \left\|
    \left( -A \right)^{ (\gamma-\theta) }
    e^{ \frac{A}{2}(t-\m{s}) }
  \right\|_{ L(H) }^2
  \bigg\|
    F'( \XY{\fl{s}} )
    \int_{ \fl{s} }^s
    \left(
      A \XY{\fl{s}}
      +
      F( \XY{\fl{s}} )
    \right) du
  \bigg\|_{ L^p( \Omega; H ) }^2
\\&+
  \K
  \left\|
    \left( -A \right)^{ (\gamma-\theta) }
    e^{ \frac{A}{2}(t-\m{s}) }
  \right\|_{ L(H) }^2
  \bigg\|
    F'( \XY{\fl{s}} )
    \int_{ \fl{s} }^s
    B( \XY{\fl{s}} ) \, dW_u
  \bigg\|_{ L^p( \Omega; H ) }^2
\\&+
  \K
  \left\|
    \left( -A \right)^{ (\gamma-\theta) }
    e^{ \frac{A}{2}(t-\m{s}) }
  \right\|_{ L(H) }^2
  \bigg[
    \int_{ \fl{s} }^s \!
    \bigg\|
      \left\|
        F''( \XY{\fl{s}} )
      \right\|_{ L^{(2)}( H, H ) }
      \left\|
        B( \XY{\fl{s}} )
      \right\|_{ HS( U_0, H ) }^2
    \bigg\|_{ L^p( \Omega; \mathbb{R} ) } \!\! du
  \bigg]^{2}
\end{split}
\end{align}
and Lemma~7.7 in Da Prato and Zabczyk \cite{dz92} hence shows
\begin{align}
\label{eq:drift_pre}
\begin{split}
 &\left\| 
    e^{ \frac{A}{2}(t-\m{s}) }
    \TF{s}
  \right\|_{ L^p( \Omega; H_{(\gamma-\theta)} ) }^2
\\&\leq
  \K
  \left(
    1
    +
    \left\|
      \XY{\fl{s}}
    \right\|_{ L^p( \Omega; H_{\alpha} ) }^2
  \right)
  +
  \frac{ \K \left( s - \fl{s} \right)}{ \left( t - \m{s} \right)^{ 2(\gamma-\theta) } }
  \int_{ \fl{s} }^s
  \left\|
    A \XY{\fl{s}}
    +
    F( \XY{\fl{s}} )
  \right\|_{ L^p( \Omega; H ) }^2 du
\\&+
  \frac{ \K }{ \left( t - \m{s} \right)^{ 2(\gamma-\theta) } }
  \int_{ \fl{s} }^s
  \left\|
    B( \XY{\fl{s}} )
  \right\|_{ L^p( \Omega; H ) }^2 du
  +
  \frac{ \K }{ \left( t - \m{s} \right)^{ 2(\gamma-\theta) } }
  \bigg[
    \int_{ \fl{s} }^s
    \bigg(
      1
      +
      \left\|
        \XY{\fl{s}}
      \right\|_{ L^p( \Omega; H ) }
    \bigg) \, du
  \bigg]^{2}
\\&\leq
  \K
  \left(
    1
    +
    \left\|
      \XY{\fl{s}}
    \right\|_{ L^p( \Omega; H_{\gamma} ) }^2
  \right)
  +
  \frac{ \K \left( s - \fl{s} \right)}{ \left( t - \m{s} \right)^{ 2(\gamma-\theta) } }
  \left(
    1
    +
    \left\|
      \XY{\fl{s}}
    \right\|_{ L^p( \Omega; H_{\gamma} ) }^2
  \right)
\\&\leq
  \K
  \left(
    1
    +
    \left\|
      \XY{\fl{s}}
    \right\|_{ L^p( \Omega; H_{\gamma} ) }^2
  \right)
\end{split}
\end{align}
for all $ s \in [0,t) $ 
and all 
$ t \in \{ 0, \tfrac{T}{M}, \ldots ,\tfrac{(M-1)T}{M}, T \} $. 
In the next step we combine
the estimate
\begin{align}
\label{eq:discrete_integral}
\begin{split}
&
  \int_{ 0 }^{ t }
  \left(
    t - \m{s}
  \right)^{ -\theta } ds
  =
  \left[
    \frac{ T }{ M }
  \right]^{ - \theta }
  \sum_{ l=0 }^{ \frac{tM}{T}-1 }
  \int_{ \frac{lT}{M} }^{ \frac{(l+1)T}{M} }
  \left(
    \tfrac{tM}{T} - l -\tfrac{1}{2}
  \right)^{ -\theta } ds
\\&=
  \left[
    \frac{ T }{ M }
  \right]^{ ( 1 - \theta ) }
  \sum_{ l=0 }^{ \frac{tM}{T}-1 }
  \left(
    l + \tfrac{1}{2}
  \right)^{ -\theta }
\leq
  \left[
    \frac{ T }{ M }
  \right]^{ ( 1 - \theta ) }
  \sum_{ l=0 }^{ \frac{tM}{T}-1 }
  \left(
    \tfrac{ l + 1}{ 2 }
  \right)^{ -\theta }
\leq
  2
  \left[
    \frac{ T }{ M }
  \right]^{ ( 1 - \theta ) }
  \sum_{ l = 1 }^{ \frac{ t M}{ T } }
  l^{ -\theta }
\\ & \leq
  \frac{\K}{ M^{ (1-\theta) }}
  \left(
    1
    +
    \int_1^M
    s^{-\theta} \, ds
  \right)
=
  \frac{\K}{ M^{ (1-\theta) } }
  \left(
    1
    +
    \frac{ ( M^{(1-\theta)} - 1 ) }{ (1-\theta) }
  \right)
  \leq
  \K
\end{split}
\end{align}
for all 
$ 
  t \in \{ 0, \tfrac{T}{M}, \ldots ,\tfrac{(M-1)T}{M}, T \} 
$ 
with inequality~\eqref{eq:drift_pre} 
to obtain that
\begin{align}
\label{eq:bounded_drift}
\begin{split}
 &\left\| 
    \int_{ 0 }^{ t }
    e^{ A(t-\m{s}) }
    \TF{s} \, ds
  \right\|_{ L^p( \Omega; H_{\gamma} ) }^2
\\&\leq 
  \left[
    \int_{ 0 }^{ t }
    2^{\theta}
    \left(
      t - \m{s}
    \right)^{ -\theta }
    \left\| 
      e^{ \frac{A}{2}(t-\m{s}) }
      \TF{s}
    \right\|_{ L^p( \Omega; H_{(\gamma-\theta)} ) } ds
  \right]^2
\\&\leq 
  \left[
    \int_{ 0 }^{ t }
    2^{\theta}
    \left(
      t - \m{s}
    \right)^{ -\theta } ds
  \right]
  \left[
    \int_{ 0 }^{ t }
    2^{\theta}
    \left(
      t - \m{s}
    \right)^{ -\theta } 
    \left\| 
      e^{ \frac{A}{2}(t-\m{s}) }
      \TF{s}
    \right\|_{ L^p( \Omega; H_{(\gamma-\theta)} ) }^2 ds
  \right]
\\&\leq
  \K
  \int_{ 0 }^{ t }
  \left(
    t - \m{s}
  \right)^{ -\theta } 
  \left(
  1
  +
  \left\|
    \XY{\fl{s}}
  \right\|_{ L^p( \Omega; H_{\gamma} ) }^2
  \right) ds
\end{split}
\end{align}
for all $ s \in [0,t) $ and all 
$ t \in \{ 0, \tfrac{T}{M}, \ldots ,\tfrac{(M-1)T}{M}, T \} 
$. 
In addition, note that
\begin{align}
\notag
 &\left\|
     e^{ \frac{A}{2}(t-\m{s}) }
     \left(
       \TB{s}
       +
       \left(
         \m{s} - s
       \right) 
       A B( \XY{\fl{s}} )
     \right)
   \right\|_{ L^p( \Omega; HS( U_0, H_{ (\gamma-\theta/2) } ) ) }^2
\\&\leq
\notag
  \K
  \left\|
    \left( -A \right)^{ (\gamma-\frac{\theta}{2}-\beta) }
  \right\|_{ L(H) }^2
  \left\|
    B( \XY{\fl{s}} )
  \right\|_{ L^p( \Omega; HS( U_0, H_{\beta} ) ) }^2
\\&+
\notag
  \K
  \left( \m{s} - s \right)^2
  \left\|
    \left( -A \right)^{ (\gamma-\frac{\theta}{2}-\beta) }
  \right\|_{ L(H) }^2
  \left\|
    A
    e^{ \frac{A}{2}(t-\m{s}) }   
  \right\|_{ L(H) }^2
  \left\|
    B( \XY{\fl{s}} )
  \right\|_{ L^p( \Omega; HS( U_0, H_{\beta} ) ) }^2
\\&+
  \K
  \left\|
    \left( -A \right)^{ (\gamma-\frac{\theta}{2}) }
    e^{ \frac{A}{2}(t-\m{s}) }   
  \right\|_{ L(H) }^2
  \bigg\|
    B'( \XY{\fl{s}} )
    \int_{ \fl{s} }^s
    \left(
      A \XY{\fl{s}}
      +
      F( \XY{\fl{s}} )
    \right) du
  \bigg\|_{ L^p( \Omega; HS( U_0, H ) ) }^2
\\&+
\notag
  \K
  \left\|
    \left( -A \right)^{ (\gamma-\frac{\theta}{2}-\delta) }
    e^{ \frac{A}{2}(t-\m{s}) }   
  \right\|_{ L(H) }^2
  \bigg\|
    B'( \XY{\fl{s}} )
    \int_{ \fl{s} }^s
    B( \XY{\fl{s}} ) \, dW_u
  \bigg\|_{ L^p( \Omega; HS( U_0, H_{ \delta } ) ) }^2
\\&+
\notag
  \K
  \left\|
    \left( -A \right)^{ (\gamma-\frac{\theta}{2}) }
    e^{ \frac{A}{2}(t-\m{s}) }   
  \right\|_{ L(H) }^2
  \bigg\|
    B'( \XY{\fl{s}} )
    \int_{ \fl{s} }^s
    B'( \XY{\fl{s}} )
    \int_{ \fl{s} }^u
    B( \XY{\fl{s}} ) \, dW_v \, dW_u
  \bigg\|_{ L^p( \Omega; HS( U_0, H ) ) }^2
\\&+
\notag
  \K
  \left\|
    \left( -A \right)^{ (\gamma-\frac{\theta}{2}) }
    e^{ \frac{A}{2}(t-\m{s}) }   
  \right\|_{ L(H) }^2
  \bigg\|
    B''( \XY{\fl{s}} )
    \bigg(
      \int_{ \fl{s} }^s \!\!\!
      B( \XY{\fl{s}} ) \, dW_u,
      \int_{ \fl{s} }^s \!\!\!
      B( \XY{\fl{s}} ) \, dW_u
    \bigg)
  \bigg\|_{ L^p( \Omega; HS( U_0, H ) ) }^2
\end{align}
and 
Lemma~7.7 in 
Da Prato \& Zabczyk~\cite{dz92}
hence implies
\begin{align}
\begin{split}
 &\left\|
     e^{ \frac{A}{2}(t-\m{s}) }
     \left(
       \TB{s}
       +
       \left(
         \m{s} - s
       \right) 
       A B( \XY{\fl{s}} )
     \right)
   \right\|_{ L^p( \Omega; HS( U_0, H_{ (\gamma-\theta/2) } ) ) }^2
\\&\leq
  \K
  \left(
    1
    +
    \frac{ \left( \m{s} - s \right)^2 }{ \left( t - \m{s} \right)^2 }  
  \right)
  \left( 
    1
    +
    \left\|
      \XY{\fl{s}}
    \right\|_{ L^p( \Omega; H_{\beta} ) }^2
  \right)
\\&+
  \K
  \frac{ \left( \fl{s} - s \right) }{ \left( t - \m{s} \right)^{ 2(\gamma-\frac{\theta}{2}) } }  
  \int_{ \fl{s} }^s
  \left\|
    A \XY{\fl{s}}
    +
    F( \XY{\fl{s}} )
  \right\|_{ L^p( \Omega; H ) }^2 du
\\&+
  \K
  \left( t - \m{s} \right)^{ -2(\gamma-\frac{\theta}{2}-\delta) }
  \int_{ \fl{s} }^s
  \left\|
    B( \XY{\fl{s}} ) 
  \right\|_{ L^p( \Omega; HS( U_0, H_{ \delta } ) ) }^2 du
\\&+
  \K
  \left( t - \m{s} \right)^{ -2(\gamma-\frac{\theta}{2}) }
  \int_{ \fl{s} }^s
  \int_{ \fl{s} }^u
  \left\|
    B( \XY{\fl{s}} )
  \right\|_{ L^p( \Omega; HS( U_0, H ) ) }^2 dv \, du
\\&+
  \K
  \left( t - \m{s} \right)^{ -2(\gamma-\frac{\theta}{2}) }
  \bigg\|
    \left\|
      B''( \XY{\fl{s}} )
    \right\|_{ L^{(2)}( H, HS(U_0,H) ) }
    \Big\|
      \int_{ \fl{s} }^s \!\!\!
      B( \XY{\fl{s}} ) \, dW_u
    \Big\|_H^2
  \bigg\|_{ L^p( \Omega; \mathbb{R} ) }^2
\end{split}
\end{align}
for all $ s \in [0,t) $ and all $ t \in \{ 0, \tfrac{T}{M}, \ldots ,\tfrac{(M-1)T}{M}, T \} $. 
This and 
again 
Lemma~7.7 in 
Da Prato \& 
Zabczyk~\cite{dz92} 
imply
\begin{align}
\label{eq:bounded_diff_pre}
\begin{split}
 &\left\|
     e^{ \frac{A}{2}(t-\m{s}) }
     \left(
       \TB{s}
       +
       \left(
         \m{s} - s
       \right) 
       A B( \XY{\fl{s}} )
     \right)
   \right\|_{ L^p( \Omega; HS( U_0, H_{ (\gamma-\theta/2) } ) ) }^2
\\&\leq
  \K
  \left( 
    1
    +
    \left\|
      \XY{\fl{s}}
    \right\|_{ L^p( \Omega; H_{\gamma} ) }^2
  \right)
\\&+
  \K
  \left( t - \m{s} \right)^{ -2(\gamma-\frac{\theta}{2}) }
  \bigg\|
    \int_{ \fl{s} }^s \!
    \left\|
      B''( \XY{\fl{s}} )
    \right\|_{ L^{(2)}( H, HS(U_0,H) ) }^{ 1/2 }
    B( \XY{\fl{s}} ) \, dW_u
  \bigg\|_{ L^{2p}( \Omega; H ) }^4
\\&\leq
  \K
  \left( 
    1
    +
    \left\|
      \XY{\fl{s}}
    \right\|_{ L^p( \Omega; H_{\gamma} ) }^2
  \right)
\\&+
  \K
  \left( t - \m{s} \right)^{ -2(\gamma-\frac{\theta}{2}) }
  \bigg[
    \int_{ \fl{s} }^s
    \bigg\|
      \left\|
        B''( \XY{\fl{s}} )
      \right\|_{ L^{(2)}( H, HS(U_0,H) ) }
      \left\|
        B( \XY{\fl{s}} )
      \right\|_{ HS( U_0, H ) }^2
    \bigg\|_{ L^{p}( \Omega; \mathbb{R} ) } du
  \bigg]^2
\\&\leq
  \K
  \left( 
    1
    +
    \left\|
      \XY{\fl{s}}
    \right\|_{ L^p( \Omega; H_{\gamma} ) }^2
  \right)
  +
  \frac{ \K \left( \fl{s} - s \right)^2 }{ \left( t - \m{s} \right)^{ 2(\gamma-\frac{\theta}{2}) } }  
  \left( 
    1
    +
    \left\|
      \XY{\fl{s}}
    \right\|_{ L^p( \Omega; H ) }^2
  \right)
\\&\leq
  \K
  \left( 
    1
    +
    \left\|
      \XY{\fl{s}}
    \right\|_{ L^p( \Omega; H_{\gamma} ) }^2
  \right)
\end{split}
\end{align}
for all 
$ s \in [0,t) $ and all 
$ t \in \{ 0, \tfrac{T}{M}, \ldots ,\tfrac{(M-1)T}{M}, T \} 
$. 
Again 
Lemma~7.7 in Da Prato 
and Zabczyk~\cite{dz92} 
hence shows 
\begin{align}
\label{eq:bounded_diff}
\begin{split}
 &\left\| 
    \int_{ 0 }^{ t }
    e^{ A(t-\m{s}) }
    \left(
       \TB{s}
       +
       \left(
         \m{s} - s
       \right) 
       A B( \XY{\fl{s}} )
     \right) dW_s
  \right\|_{ L^p( \Omega; H_{\gamma} ) }^2
\\&\leq
  \K
  \int_{ 0 }^{ t }
  \left(
    t - \m{s}
  \right)^{ -\theta }
  \left\|
     e^{ \frac{A}{2}(t-\m{s}) }
     \left(
       \TB{s}
       +
       \left(
         \m{s} - s
       \right) 
       A B( \XY{\fl{s}} )
     \right)
   \right\|_{ L^p( \Omega; HS( U_0, H_{ (\gamma-\theta/2) } ) ) }^2 ds
\\&\leq
  \K
  \int_{ 0 }^{ t }
  \left(
    t - \m{s}
  \right)^{ -\theta }
  \left( 
    1
    +
    \left\|
      \XY{\fl{s}}
    \right\|_{ L^p( \Omega; H_{\gamma} ) }^2
  \right) ds
\end{split}
\end{align}
for all $ s \in [0,t) $ and all 
$ t \in \{ 0, \tfrac{T}{M}, \ldots ,\tfrac{(M-1)T}{M}, T \} 
$. 
Next we combine 
\eqref{eq:scheme_new}, 
\eqref{eq:bounded_drift}, 
\eqref{eq:bounded_diff} 
and \eqref{eq:discrete_integral}
to obtain that
\begin{align}
\label{eq:bounded_main}
\begin{split}
&
  \left\| 
    \Y{m} 
  \right\|_{ L^p( \Omega; H_{\gamma} ) }^2
\leq
  \K
  \left\| 
    e^{ A\frac{mT}{M} } \xi
  \right\|_{ L^p( \Omega; H_{\gamma} ) }^2
  +
  \K
  \left\| 
    \int_{ 0 }^{ \frac{mT}{M} }
    e^{ A(\frac{mT}{M}-\m{s}) }
    \TF{s} \, ds
  \right\|_{ L^p( \Omega; H_{\gamma} ) }^2
\\&\quad+
  \K
  \left\| 
    \int_{ 0 }^{ \frac{mT}{M} }
    e^{ A(\frac{mT}{M}-\m{s}) }
    \left(
       \TB{s}
       +
       \left(
         \m{s} - s
       \right) 
       A B( \XY{\fl{s}} )
     \right) dW_s
  \right\|_{ L^p( \Omega; H_{\gamma} ) }^2
\\&\quad\leq
  \K
  \left\| 
    \xi
  \right\|_{ L^p( \Omega; H_{\gamma} ) }^2
  +
  \K
  \int_{ 0 }^{ \frac{mT}{M} }
  \left(
    \tfrac{mT}{M} - \m{s}
  \right)^{ -\theta }
  \left( 
    1
    +
    \left\|
      \XY{\fl{s}}
    \right\|_{ L^p( \Omega; H_{\gamma} ) }^2
  \right) ds
\\ & \quad \leq
  \K
  \left(
    1 +
    \left\| 
      \xi
    \right\|_{ L^p( \Omega; H_{\gamma} ) }^2 
  \right)
  +
  \K
  \int_{ 0 }^{ \frac{mT}{M} }
  \frac{
    2^{ \theta }
  }{
    \left(
      2 \left[
        \tfrac{mT}{M} - \m{s}
      \right]
    \right)^{ \theta }
  }
    \left\|
      \XY{\fl{s}}
    \right\|_{ L^p( \Omega; H_{\gamma} ) }^2
  ds
\\& \quad =
  \K
  \left(
    1 +
    \left\| 
      \xi
    \right\|_{ L^p( \Omega; H_{\gamma} ) }^2 
  \right)
  +
  \K 
  \left[ 
    \frac{ T }{ M }
  \right]^{
    ( 1 - \theta )
  }
  \sum_{ l=0 }^{ m-1 }
  \frac{
    2^{ \theta }
  }{
  \left(
    m - l + 
    \left[ m - l - 1 \right]
  \right)^{ \theta }
  }
    \left\|
      \Y{l}
    \right\|_{ L^p( \Omega; H_{\gamma} ) }^2
\\&\quad\leq
  \K
  \left(
    1
    +
    \left\| 
      \xi
    \right\|_{ L^p( \Omega; H_{\gamma} ) }^2
  \right)
  +
  \K M^{ -(1-\theta) }
  \sum_{ l=0 }^{ m-1 }
  \left( m - l \right)^{ -\theta }
  \left\|
    \Y{l}
  \right\|_{ L^p( \Omega; H_{\gamma} ) }^2
\end{split}
\end{align}
for all $ m \in \{ 0, 1, \ldots, M \} $ 
and all $ M \in \mathbb{N} $. 
In the next step we use the mappings 
$ 
  E_{ \varepsilon } 
  \colon [0, \infty) \rightarrow [0, \infty) 
$,
$ 
  \varepsilon \in ( 0, \infty )
$,
defined by 
$ 
  E_{ \varepsilon} (t) 
:= 
  \sum_{ n=0 }^{ \infty } 
  \frac{ t^{n\varepsilon} }{ \Gamma(n\varepsilon + 1) }
$ 
for all 
$ \varepsilon \in ( 0, \infty) $ 
and all $ t \in [0, \infty) $ 
(see Section~7 in Henry \cite{h81}) and 
apply a generalized version of the discrete Gronwall lemma (see Theorem~6.1 in Dixon \& 
McKee~\cite{dm86}) to \eqref{eq:bounded_main} 
to obtain
\begin{align}
\label{eq:w_1}
\begin{split}
  \left\| 
    \Y{m} 
  \right\|_{ L^p( \Omega; H_{\gamma} ) }^2
 &\leq
  \K
  \left(
    1
    +
    \left\| 
      \xi
    \right\|_{ L^p( \Omega; H_{\gamma} ) }^2
  \right)
  E_{ (1-\theta) }\!\left(2 M \left(\K M^{ -(1-\theta) }\Gamma(1-\theta)\right)^{\frac{1}{(1-\theta)} }  \right)
\\&=
  \K
  \left(
    1
    +
    \left\| 
      \xi
    \right\|_{ L^p( \Omega; H_{\gamma} ) }^2
  \right)
  E_{ (1-\theta) }\!\left(2 \left(\K \Gamma(1-\theta)\right)^{\frac{1}{(1-\theta)} }  \right)
\\&\leq
  \K
  \left(
    1
    +
    \left\| 
      \xi
    \right\|_{ L^p( \Omega; H_{\gamma} ) }^2
  \right)
  E_{ (1-\theta) }\!\left(
    2 \left( \K \Gamma( \tfrac{ 1 }{ K } ) \right)^K
  \right)
\\ & \leq 
  \K
  \left(
    1
    +
    \left\| 
      \xi
    \right\|_{ L^p( \Omega; H_{\gamma} ) }^2
  \right)
    E_{ (1-\theta) }\!\left(
      \K 
  \right)
\end{split}
\end{align}
for all $ m \in \{ 0, 1, \ldots, M \} $ 
and all $ M \in \mathbb{N} $.
In addition, note that
\begin{align}
\label{eq:w_2}
\begin{split}
  E_{ (1-\theta) }\!\left(
      \K 
  \right)
 &=
  \sum_{ n=0 }^{ \infty } 
  \frac{ \K^{n(1-\theta)} }{ \Gamma(n(1-\theta) + 1) }
  =
  \sum_{ n=0 }^{ \cll{\frac{1}{K}} } 
  \frac{ \K^{n(1-\theta)} }{ \Gamma(n(1-\theta) + 1) }
  +
  \sum_{ n=\cll{\frac{1}{K}}+1 }^{ \infty } 
  \frac{ \K^{n(1-\theta)} }{ \Gamma(n(1-\theta) + 1) }
\\&\leq
  \left(
    \cll{\frac{1}{K}} + 1
  \right)
  \K^{ ( 2 + \frac{1}{K} ) }
  +
  \sum_{ n=\cll{\frac{1}{K}}+1 }^{ \infty } 
  \frac{ \K^{\frac{n}{K}} }{ \Gamma(\frac{n}{K} + 1) }
  \leq
  \K
  +
  \sum_{ n=0 }^{ \infty } 
  \frac{ \K^{\frac{n}{K}} }{ \Gamma(\frac{n}{K} + 1) }
\\&=
  \K
  +
  E_{ \frac{1}{K} }\!\left(
      \K 
  \right)
  \leq
  \K.
\end{split}
\end{align}
Combining \eqref{eq:w_1} and \eqref{eq:w_2} then gives
\begin{align}
\begin{split}
  \left\| 
    \Y{m} 
  \right\|_{ L^p( \Omega; H_{\gamma} ) }^2
\leq
  \K
  \left(
    1
    +
    \left\| 
      \xi
    \right\|_{ L^p( \Omega; H_{\gamma} ) }^2
  \right)
\end{split}
\end{align}
for all $ m \in \{ 0, 1, \ldots, M \} $ 
and all $ M \in \mathbb{N} $.
This completes the proof of Proposition \ref{prop:well_defined}.

\subsection{Proof of Theorem~\ref{thm:main}}
\label{sec:proof_of_thm1}
Throughout this proof 
$ C \colon [0,\infty) \to [0,\infty) $
is a universal non-decreasing
function which changes from line to line.
Note that Jensen's inequality implies
\begin{equation}
\label{eq:boundedness}
  \left\|
    \XY{t}
  \right\|_{ L^p( \Omega; H ) }^2
  \leq
  \left\|
    \XY{t}
  \right\|_{ L^6( \Omega; H ) }^2
  \leq
  C(K)
  \left(
    1
    +
    \left\|
      \XY{\fl{t}}
    \right\|_{ L^6( \Omega; H ) }^2
  \right)
  \leq
  C(K)
  \left(
    1
    +
    \left\|
      \XY{\fl{t}}
    \right\|_{ L^6( \Omega; H_{ \gamma } ) }^2
  \right)
  \leq
  C(K)
\end{equation}
and
\begin{equation}
\label{eq:hoelder}
  \left\|
    \XY{t} - \XY{\fl{t}}
  \right\|_{ L^p( \Omega; H ) }
  \leq
  \left\|
    \XY{t} - \XY{\fl{t}}
  \right\|_{ L^6( \Omega; H ) }
  \leq
  \C
  \left(t-\fl{t} \right)^{ 1/2 }
\end{equation}
for all $ p \in (0,6] $, $ t \in [0, T] $ and all $ M \in \mathbb{N} $. Moreover, let $ \Z{} \colon [0,T] \times \Omega \rightarrow H $, $ M \in \mathbb{N} $, be stochastic processes satisfying
\begin{align}
\begin{split}
\label{eq:helper_process_Z}
  \Z{t}
 &=
  e^{ At } \X{0}
  +
  \int_{ 0 }^{ t }
  e^{ A(t - s) }
  F( \XY{s} ) \, ds
  +
  \int_{ 0 }^{ t }
  e^{ A(t - s) }
  B( \XY{s} ) \, dW_s
\end{split}
\end{align}
$ \mathbb{P} $-a.s.\ for all $ t \in [0,T] $ and all $ M \in \mathbb{N} $ and let $ \TZ{} \colon \{ 0,1,\ldots,M \} \times \Omega \rightarrow H $, $ M \in \mathbb{N} $, be stochastic processes satisfying
\begin{align}
\begin{split}
\label{eq:helper_Z_2}
  \TZ{t}
 &=
  e^{ At } X_0
  +
  \int_{ 0 }^{ t }
  e^{ A(t - s) } \,
  \TF{s} \, ds
  +
  \int_{ 0 }^{ t }
  e^{ A(t - s) } \,
  \TB{s} \, dW_s
\end{split}
\end{align}
$ \mathbb{P} $-a.s.\ for all $ t \in [0,T] $ and all $ M \in \mathbb{N} $. Next observe that the triangle inequality implies that
\begin{align}
\begin{split}
\label{eq:start}
  \Big\|
    \X{ \frac{mT}{M} } - \Y{ m }
  \Big\|_{ L^2( \Omega; H ) }^2
 &\leq
  2
  \Big\|
    \X{ \frac{mT}{M} } - \Z{ \frac{mT}{M} }
  \Big\|_{ L^2( \Omega; H ) }^2
  +
  2
  \left[
    \left\|
      \Z{ \frac{mT}{M} } - \TZ{ \frac{mT}{M} }
    \right\|_{ L^2( \Omega; H ) }
    +
    \left\|
      \TZ{ \frac{mT}{M} } - \Y{ m }
    \right\|_{ L^2( \Omega; H ) }
  \right]^2
\end{split}
\end{align}
for all $ m \in \{ 0,1, \ldots, M \} $ and all $ M \in \mathbb{N} $. Moreover, note that
\begin{align}
\notag
  \left\|
    \X{t}
    -
    \Z{t}
  \right\|_{ L^2( \Omega; H ) }^2
 &\leq
  \C
  \left\|
    \int_{ 0 }^{ t }
    e^{ A( t-s ) }
    \left( F( \X{s} ) - F( \XY{s} ) \right) ds
  \right\|_{ L^2( \Omega; H ) }^2
\\&+
\label{eq:lipschitz}
  \C
  \left\|
    \int_{ 0 }^{ t }
    e^{ A( t-s ) }
    \left( B( \X{s} ) - B( \XY{s} ) \right) dW_s
  \right\|_{ L^2( \Omega; H ) }^2
\\&\leq
\notag
  \C
  \int_{ 0 }^{ t }
  \left( 
  \left\|
     F( \X{s} ) - F( \XY{s} )
  \right\|_{ L^2( \Omega; H ) }^2
  +
  \left\|
    B( \X{s} ) - B( \XY{s} )
  \right\|_{ L^2( \Omega; HS( U_0, H ) ) }^2
  \right) ds
\\&\leq
\notag
  \C
  \int_{ 0 }^{ t }
  \left\|
     \X{s} - \XY{s}
  \right\|_{ L^2( \Omega; H ) }^2 ds 
\end{align}
for all $ t \in [0,T] $ and all $ M \in \mathbb{N} $. In addition, \eqref{eq:helper_process_Z}, \eqref{eq:helper_Z_2}, \eqref{eq:tay_F_end} and \eqref{eq:tay_B_end} imply
\begin{align}
\label{eq:est_taylor_F}
  \left\|
    \Z{ \frac{mT}{M} } - \TZ{ \frac{mT}{M}  }
  \right\|_{ L^2( \Omega; H ) }
 &\leq
  \left\|
    \int_{ 0 }^{ \frac{mT}{M}  }
    e^{A(\frac{mT}{M} - s)}
    \left(
      F( \XY{s} )
      -
      \TF{s}
    \right) ds
  \right\|_{ L^2( \Omega; H ) }
\\&\quad+
\label{eq:est_taylor_B}
  \left[
    \int_{ 0 }^{ \frac{mT}{M}  }
    \left\|
      B( \XY{s} )
      -
      \TB{s}
    \right\|_{ L^2( \Omega; HS( U_0, H ) ) }^2 ds
  \right]^{ \frac{1}{2} }
\\&\leq
  \C M^{ -\gamma }
  +
  \left[
    \int_{ 0 }^{ \frac{mT}{M}  }
    \C M^{ -2\gamma } \, ds
  \right]^{ \frac{1}{2} }
  \leq
\label{eq:tt_end}
  \C M^{ -\gamma }
\end{align}
for all $ m \in \{ 0,1, \ldots, M \} $ and all $ M \in \mathbb{N} $. Combining \eqref{eq:start}--\eqref{eq:tt_end} yields 
\begin{align}
\label{eq:mid_est}
\begin{split}
  \Big\|
    \X{ \frac{mT}{M} } - \Y{ m }
  \Big\|_{ L^2( \Omega; H ) }^2
 &\leq
  \C
  \int_{ 0 }^{ \frac{mT}{M} }
  \left\|
     \X{s} - \XY{s}
  \right\|_{ L^2( \Omega; H ) }^2 ds 
  +
  2
  \left[
    \C M^{ -\gamma }
    +
    \left\|
      \TZ{ \frac{mT}{M} } - \Y{ m }
    \right\|_{ L^2( \Omega; H ) }
  \right]^2
\end{split}
\end{align}
for all $ m \in \{ 0,1, \ldots, M \} $ and all $ M \in \mathbb{N} $. Next \eqref{eq:scheme_new}, \eqref{eq:helper_Z_2}, \eqref{eq:discrete_F_end} and \eqref{eq:discrete_B_end} show
\begin{align}
 &\left\|
    \TZ{ \frac{mT}{M} } - \Y{ m }
  \right\|_{ L^2( \Omega; H ) }
\\&\leq
\label{eq:discrete_F}
  \left\|
    \int_{ 0 }^{ \frac{mT}{M} }
    \left( 
      e^{A(\frac{mT}{M}-s)}
      -
      e^{ A(\frac{mT}{M}-\m{s}) }
    \right) 
    \TF{s} \, ds
  \right\|_{ L^2( \Omega; H ) }
\\&+
\label{eq:discrete_B}
  \left\|
    \int_{ 0 }^{ \frac{mT}{M} }
    \left( 
      e^{A(\frac{mT}{M}-s)} \,
      \TB{s}
      -
      e^{A(\frac{mT}{M}-\m{s})}
      \left(
        \TB{s}
        +
        \left( \m{s} \!-\! s \right)\! A
        B( \XY{\fl{s}} )
      \right)
    \right) dW_s
  \right\|_{ L^2( \Omega; H ) }^2
\\&\leq
  \C M^{ -\min(\alpha+1,\frac{3}{2}) }
  \left( 1 + \log(M) \right)
  +
  \C M^{ -\min(\beta+\frac{1}{2},\delta+1,\frac{3}{2}) }
  \left( 1 + \log(M) \right)
  \leq
  \C M^{ -\gamma }
\end{align}
for all $ m \in \{ 0,1, \ldots, M \} $ and all $ M \in \mathbb{N} $. Combining this and \eqref{eq:mid_est} then yields
\begin{align}
\label{eq:main_3}
\begin{split}
  \Big\|
    \X{ \frac{mT}{M} } - \Y{ m }
  \Big\|_{ L^2( \Omega; H ) }^2
 &\leq
  \C
  \int_{ 0 }^{ \frac{mT}{M} }
  \left\|
     \X{s} - \XY{s}
  \right\|_{ L^2( \Omega; H ) }^2 ds 
  + 
  \C M^{ -2\gamma }
\\&=
  \C
  \sum_{ l=0 }^{ m-1 }
  \int_{ \frac{lT}{M} }^{ \frac{(l+1)T}{M} }
  \left\|
     \X{s} - \XY{s}
  \right\|_{ L^2( \Omega; H ) }^2 ds 
  + 
  \C M^{ -2\gamma }
\\&\leq
  \C
  \sum_{ l=0 }^{ m-1 }
  \int_{ \frac{lT}{M} }^{ \frac{(l+1)T}{M} }
  \left\|
     \X{\frac{lT}{M}} - \Y{l}
  \right\|_{ L^2( \Omega; H ) }^2 ds 
  + 
  \C M^{ -2\gamma }
\end{split}
\end{align}
for all $ m \in \{ 0,1, \ldots, M \} $ and all $ M \in \mathbb{N} $. To finish the proof of Theorem \ref{thm:main} we apply the discrete Gronwall lemma to \eqref{eq:main_3} and take square root to obtain
\begin{equation}
\label{eq:main_goal}
  \Big\|
    \X{ \frac{mT}{M} } - \Y{ m }
  \Big\|_{ L^2( \Omega; H ) }
  \leq
  \C M^{ -\gamma }
\end{equation}
for all $ m \in \{ 0,1, \ldots, M \} $ and all $ M \in \mathbb{N} $.

\subsubsection{Estimates for $ \big\| \Z{m} - \TZ{m} \big\|_{ L^2( \Omega; H ) }$ for $ m \in \{ 0,1, \ldots, M \} $ and $ M \in \mathbb{N} $}
The following well known lemma will be used frequently below.
\begin{lemma}
\label{lemma_eA}
Let the setting in Section~\ref{sec:sec2} be fulfilled. Then
\begin{equation}
 \left\|
   \left( -tA \right)^{-\kappa} 
   \left( 
     e^{ At } - I
   \right)
 \right\|_{ L(H) }
 \leq
 1
\end{equation}
for all $ t \in (0, \infty) $ and all $ \kappa \in [0,1] $ and
\begin{equation}
\label{eq:lemm_21}
 \left\|
    \left( - t A \right)^{-\kappa}
    \left(
      e^{ A t }
      -
      I
      -
      t A
    \right)
  \right\|_{L(H)}
  \leq
  1
\end{equation}
for all $ t \in (0, \infty) $ and all $ \kappa \in [1,2] $.
\end{lemma}
%
%
%

With the help of Lemma~\ref{lemma_eA} we first establish some estimates 
that we exploit in the estimation 
of \eqref{eq:est_taylor_F} and \eqref{eq:est_taylor_B}. 
More formally, observe that Lemma~\ref{lemma_eA} implies
\begin{align}
\begin{split}
 &\left\|
    \int_{ \fl{t} }^t
    \left(
      e^{ A(t-s) }
      F( \XY{s} )
      -
      F( \XY{\fl{s}} )
    \right) ds
  \right\|_{ L^2( \Omega; H ) }
\\&\leq
  \C
  \int_{ \fl{t} }^t
  \left(
    \left\|
      e^{ A( t-s ) }
      \left( 
        F( \XY{s} )
        -
        F( \XY{\fl{s}} )
      \right)
    \right\|_{ L^2( \Omega; H ) }
    +
    \left\|
      \left(
        e^{ A( t-s ) }
        -
        I
      \right)
      F( \XY{\fl{s}} )
    \right\|_{ L^2( \Omega; H ) }
  \right) ds 
\\&\leq
  \C
  \int_{ \fl{t} }^t
  \left(
    \left\|
      \XY{s}
      -
      \XY{\fl{s}}
    \right\|_{ L^2( \Omega; H ) }
    +
    \left\|
      \left( -A \right)^{ -\alpha }
      \left(
        e^{ A( t-s ) }
        -
        I
      \right)
    \right\|_{ L(H) }
    \left\|
      F( \XY{\fl{s}} )
    \right\|_{ L^2( \Omega; H_{\alpha} ) }
  \right) ds
\\&\leq
\label{eq:drift_term}
  \C
  \int_{ \fl{t} }^t
  \left(
    M^{-\frac{1}{2} }
    +
    M^{ - \min( \alpha, 1 ) }
    \left(
      1
      +
      \left\| \XY{\fl{s}} \right\|_{ L^2( \Omega; H_{\alpha} ) }
    \right)
  \right) ds
  \leq
  \C
  M^{ -\min(\alpha+1, \frac{3}{2}) }
  \leq
  \C M^{ -\gamma }  
\end{split}
\end{align}
and
\begin{align}
\notag
 &\left\|
    \int_{ \fl{t} }^t
    \left(
      e^{ A(t-s) }
      B( \XY{s} )
      -
      B( \XY{\fl{s}} )
    \right) dW_s
  \right\|_{ L^p( \Omega; H ) }^2
\\&\leq
\notag
  \C
  \int_{ \fl{t} }^t
  \left\|
    e^{ A( t-s ) }
    \left( 
      B( \XY{s} )
      -
      B( \XY{\fl{s}} )
    \right)
    +
    \left(
      e^{ A( t-s ) }
      -
      I
    \right)
    B( \XY{\fl{s}} )
  \right\|_{ L^p( \Omega; HS( U_0, H ) ) }^2 ds 
\\&\leq
\notag
  \C
  \int_{ \fl{t} }^t
  \left(
    \left\|
      \XY{s}
      -
      \XY{\fl{s}}
    \right\|_{ L^p( \Omega; H ) }^2
    +
    \left\|
      \left( -A \right)^{ -\beta }
      \left(
        e^{ A( t-s ) }
        -
        I
      \right)
    \right\|_{ L(H) }^2
    \left\|
      B( \XY{\fl{s}} )
    \right\|_{ L^p( \Omega; HS( U_0, H_{ \beta } ) ) }^2 
  \right) ds
\\&\leq
\label{eq:noise_term}
  \C
  \int_{ \fl{t} }^t
  \left(
    M^{-1}
    +
    M^{ - \min( 2 \beta , 2 ) }
    \left(
      1
      +
      \left\| \XY{\fl{s}} \right\|_{ L^p( \Omega; H_{ \beta } ) }^2
    \right)
  \right) ds
  \leq
  \C
  M^{ -\min(2\beta+1,2) }
  \leq
  \C
  M^{ -2 }
\end{align}
for all $ t \in [0, T] $, all $ M \in \mathbb{N} $ and all $ p \in [2,6] $. 
Additionally, 
Lemma~\ref{lemma_eA},
\eqref{eq:drift_term}
and \eqref{eq:hoelder} 
show
\begin{align}
\label{eq:det_solution}
\begin{split}
 &\left\|
    \XY{t}
    -
    \XY{\fl{t}}
    -
    \int_{ \fl{t} }^t
    \left(
      A \XY{\fl{s}}
      +
      F( \XY{\fl{s}} )
    \right) ds
    -
    \int_{ \fl{t} }^t
    B( \XY{s} ) \, dW_s
  \right\|_{ L^2( \Omega; H ) }
\\&\leq
  \left\|
    \left(
      e^{ A(t-\fl{t}) }
      -
      I
      -
      \left(
        t-\fl{t}
      \right) A
    \right)
    \XY{\fl{t}}
  \right\|_{ L^2( \Omega; H ) }
  +
  \left\|
    \int_{ \fl{t} }^t
    \left(
      e^{ A(t-s) }
      F( \XY{s} )
      -
      F( \XY{\fl{s}} )
    \right) ds
  \right\|_{ L^2( \Omega; H ) }
\\&+
  \left\|
    \int_{ \fl{t} }^t
    \left(
      e^{ A(t-s) }
      -
      I
    \right) 
    B( \XY{s} ) \, dW_s
  \right\|_{ L^2( \Omega; H ) }
\\&\leq
  \left\|
    \left( -A \right)^{ -\gamma }
    \left(
      e^{ A(t-\fl{t}) }
      -
      I
      -
      \left(
        t-\fl{t}
      \right) A
    \right)
  \right\|_{ L( H ) }
  \left\|
    \XY{\fl{t}}
  \right\|_{ L^2( \Omega; H_{ \gamma } ) }
  +
  \C M^{ -\gamma }
\\&+
  \left[
    \int_{ \fl{t} }^t
    \left\|
      \left( -A \right)^{ -\beta }
      \left(
        e^{ A(t-s) }
        -
        I
      \right)
    \right\|_{ L(H) }^2
    \left\|
      B( \XY{s} )
    \right\|_{ L^2( \Omega; HS( U_0, H_{ \beta } ) ) }^2 ds
  \right]^{ \frac{1}{2} }
\\&\leq
  \C M^{ -\gamma }
  +
  \C M^{ -(\beta + \frac{1}{2} ) }
  \leq
  \C M^{ -\gamma }
\end{split}
\end{align}
for all $ t \in [0, T] $ and all $ M \in \mathbb{N} $. 
Moreover, 
Lemma~\ref{lemma_eA} and \eqref{eq:noise_term} imply 
\begin{align}
\begin{split}
 &\left\|
    \XY{t}
    -
    \XY{\fl{t}}
    -
    \int_{\fl{t}}^t
    B( \XY{\fl{s}} ) \, dW_s
  \right\|_{ L^p( \Omega; H ) }
\\&\leq
  \left\|
    \left( -A \right)^{ -1 }
    \left(
      e^{ A(t-\fl{t}) }
      -
      I
    \right)
  \right\|_{ L(H) }
  \left\|
    \XY{\fl{t}}
  \right\|_{ L^p( \Omega; H_1 ) }
  +
  \int_{\fl{t}}^t
  \left\|
    e^{ A(t-s) }
    F( \XY{s} ) 
  \right\|_{ L^p( \Omega; H ) } ds
\\&
+
\label{eq:minus_noise}
  \left\|
    \int_{\fl{t}}^t
    \left( 
      e^{ A(t-s) }
      B( \XY{s} )
      -
      B( \XY{\fl{s}} )
    \right) dW_s
  \right\|_{ L^p( \Omega; H ) }
\leq
  \C M^{ -1 }
\end{split}
\end{align}
for all $ t \in [0, T] $, all $ M \in \mathbb{N} $ and all $ p \in [2,6] $. 

\paragraph{Estimation of \eqref{eq:est_taylor_F}}
\label{sec:est_taylor_F}
Note that
\begin{align}
\label{eq:taylor_exp_F}
\begin{split}
  F( \XY{t} )
 &=
  F( \XY{\fl{t}} )
  +
  F'( \XY{\fl{t}} )
  \left(
    \XY{t} - \XY{\fl{t}}
  \right)
  +
  \frac{1}{2}
  F''( \XY{\fl{t}} )
  \left(
    \XY{t} - \XY{\fl{t}},
    \XY{t} - \XY{\fl{t}}
  \right)
\\&+
  \int_0^1
  \left(
    F''\!\left( \XY{\fl{t}} + r(\XY{t} - \XY{\fl{t}}) \right)
    -
    F''( \XY{\fl{t}} )
  \right)
  \left(
    \XY{t} - \XY{\fl{t}}, \XY{t} - \XY{\fl{t}}
  \right)
  \left( 1 - r \right) dr
\end{split}
\end{align}
and
\begin{align}
\label{eq:taylor_F3}
\begin{split}
 &\left\|
    \int_0^1
    \left(
      F''\!\left( \XY{\fl{t}} + r(\XY{t} - \XY{\fl{t}}) \right)
      -
      F''( \XY{\fl{t}} )
    \right)
    \left(
      \XY{t} - \XY{\fl{t}}, \XY{t} - \XY{\fl{t}}
    \right)
    \left( 1 - r \right) dr
  \right\|_{ L^2( \Omega; H ) }
\\&\leq
  \int_0^1
  \left\|
    \left\|
      F''\!\left( \XY{\fl{t}} + r(\XY{t} - \XY{\fl{t}}) \right)
      -
      F''( \XY{\fl{t}} )
    \right\|_{ L^{(2)}(H,H) }
    \left\|
      \XY{t} - \XY{\fl{t}}
    \right\|_{ H }^2
  \right\|_{ L^2( \Omega; \mathbb{R} ) } dr
\\&\leq
  \C
  \left\|
    \XY{t} - \XY{\fl{t}}
  \right\|_{ L^6( \Omega; H ) }^3
  \leq
  \C
  M^{ -\frac{3}{2} }
\end{split}
\end{align}
for all $ t \in [0, T] $ and all $ M \in \mathbb{N} $. The remainder terms in \eqref{eq:taylor_F3} and \eqref{eq:b''_help2} are here estimated similarly as in Kruse~\cite{rk13}.
Combining \eqref{eq:taylor_exp_F} and \eqref{eq:taylor_F3} then shows
\begin{align}
\notag
 &\left\|
    \int_{ 0 }^{ t }
    e^{A(t-s)}
    \left(
      F( \XY{s} )
      -
      \TF{s}
    \right) ds
  \right\|_{ L^2( \Omega; H ) }
\\&\leq
\notag
  \left\|
    \int_{ 0 }^{ t }
    e^{A(t-s)}
    F'( \XY{\fl{s}} ) \!
    \left( \!
      \XY{s} \!-\! \XY{\fl{s}}
      \!-\!
      \int_{\fl{s}}^s
      \left(
        A \XY{\fl{u}}
        \!+\!
        F( \XY{\fl{u}} )
      \right) du
      \!-\!
      \int_{\fl{s}}^s \!
      B( \XY{\fl{u}} ) \, dW_u \!
    \right) ds
  \right\|_{ L^2( \Omega; H ) }
\\&+
\notag
  \left\|
    \int_{ 0 }^{ t }
    e^{A(t-s)}
    F''( \XY{\fl{s}} ) \!
    \left( \!
      \XY{s} \!-\! \XY{\fl{s}}
      \!-\!
      \int_{\fl{s}}^s
      B( \XY{\fl{u}} ) \, dW_u, 
      \XY{s} \!-\! \XY{\fl{s}}
      \!+\!
      \int_{\fl{s}}^s
      B( \XY{\fl{u}} ) \, dW_u \!
    \right) ds
  \right\|_{ L^2( \Omega; H ) }
\\&+
\label{eq:tay_F}
  \bigg\|
    \int_{ 0 }^{ t }
    e^{A(t-s)}
    \bigg( \!
      F''( \XY{\fl{s}} )
      \left(
        \int_{\fl{s}}^s
        B( \XY{\fl{u}} ) \, dW_u,
        \int_{\fl{s}}^s
        B( \XY{\fl{u}} ) \, dW_u
      \right)
\\&\qquad-
\notag
      \sum_{ j \in \mathcal{J} }
      \int_{\fl{s}}^s
      F''( \XY{\fl{s}} )
      \left(
        B( \XY{\fl{u}} ) g_j, B( \XY{\fl{u}} ) g_j
      \right) du
    \bigg) \, ds
  \bigg\|_{ L^2( \Omega; H ) }
  +
  \C
  M^{ -\frac{3}{2} }
\end{align}
for all $ t \in [0, T] $ and all $ M \in \mathbb{N} $. Moreover, \eqref{eq:det_solution} and \eqref{eq:hoelder} imply
\begin{align}
\notag
 &\left\|
    \int_{ 0 }^{ t }
    e^{A(t-s)}
    F'( \XY{\fl{s}} )
    \left( 
      \XY{s} - \XY{\fl{s}}
      -
      \int_{\fl{s}}^s
      \left(
        A \XY{\fl{u}}
        +
        F( \XY{\fl{u}} )
      \right) du
      -
      \int_{\fl{s}}^s
      B( \XY{\fl{u}} ) \, dW_u
    \right) ds
  \right\|_{ L^2( \Omega; H ) }
\\&\leq
\notag
  \C
  \int_{ 0 }^{ t }
  \left\|
    \XY{s}
    -
    \XY{\fl{s}}
    -
    \int_{\fl{s}}^s
    \left(
      A \XY{\fl{u}}
      +
      F( \XY{\fl{u}} )
    \right) du
    -
    \int_{\fl{s}}^s
    B( \XY{u} ) \, dW_u
  \right\|_{ L^2( \Omega; H ) } ds
\\&+
\label{eq:tay_F'}
  \C
  \left[
    \frac{T}{M}
    \sum_{ l=0 }^{ M\fl{t}/T }
    \int_{ \frac{lT}{M} }^{ \min\left(\frac{(l+1)T}{M},t\right) }
    \left\|
      \int_{ \frac{lT}{M} }^s
      \left(
        B( \XY{u} )
        -
        B( \XY{\fl{u}} )
      \right) dW_u
  \right\|_{ L^2( \Omega; H ) }^2 ds
  \right]^{ \frac{1}{2} }
\\&\leq
\notag
  \C
  \int_{ 0 }^{ t }
  M^{ -\gamma } \, ds
  +
  \C
  \left[
    M^{ -1 }
    \sum_{ l=0 }^{ M\fl{t}/T }
    \int_{ \frac{lT}{M} }^{ \min\left(\frac{(l+1)T}{M},t\right) }
    \int_{ \frac{lT}{M} }^s
    \left\|
      \XY{u}
      -
      \XY{\fl{u}}
    \right\|_{ L^2( \Omega; H ) }^2 du \, ds
  \right]^{ \frac{1}{2} }
\\&\leq
\notag
  \C M^{ -\gamma }
  +
  \C
  \left[
    \sum_{ l=0 }^{ M\fl{t}/T }
    M^{ -4 }
  \right]^{ \frac{1}{2} }
  \leq
  \C M^{ -\gamma }
\end{align}
and Lemma~7.7 in Da Prato and Zabczyk \cite{dz92}, \eqref{eq:minus_noise} and \eqref{eq:hoelder} show
\begin{align}
\notag
  &\left\|
    \int_{ 0 }^{ t }
    e^{A(t-s)}
    F''(\XY{\fl{s}})
    \left( 
      \XY{s} \!-\! \XY{\fl{s}}
      \!-\!
      \int_{\fl{s}}^s
      B( \XY{\fl{s}} ) \, dW_u, 
      \XY{s} \!-\! \XY{\fl{s}}
      \!+\!
      \int_{\fl{s}}^s
      B( \XY{\fl{u}} ) \, dW_u
    \right) ds
  \right\|_{ L^2( \Omega; H ) }
\\&\leq
\notag
  \C
  \int_{ 0 }^{ t }
  \left\|
    \XY{s}
    -
    \XY{\fl{s}}
    -
    \int_{\fl{s}}^s
    B( \XY{\fl{u}} ) \, dW_u
  \right\|_{ L^4( \Omega; H ) }
  \left\|
    \XY{s}
    -
    \XY{\fl{s}}
    +
    \int_{\fl{s}}^s
    B( \XY{\fl{u}} ) \, dW_u
  \right\|_{ L^4( \Omega; H ) } ds
\\&\leq
\label{eq:tay_F''}
  \C M^{ -1 } 
  \int_{ 0 }^{ t }
  \left(
    \bigg\|
      \XY{s}
      \!-\!
      \XY{\fl{s}}
    \bigg\|_{ L^4( \Omega; H ) }
    +
    \bigg\|
      \int_{\fl{s}}^s
      B( \XY{\fl{u}} ) \, dW_u
    \bigg\|_{ L^4( \Omega; H ) } 
  \right) ds
\\&\leq
\notag
  \C M^{ -1 } 
  \int_{ 0 }^{ t }
  \left(
    M^{ -\frac{1}{2} } 
    +
    \bigg[
      \int_{\fl{s}}^s
      \left\|
        B( \XY{\fl{u}} )
      \right\|_{ L^4( \Omega; HS( U_0, H ) }^2 du
    \bigg]^{ \frac{1}{2} }
  \right) ds
  \leq
  \C M^{ -\frac{3}{2} }
\end{align}
for all $ t \in [0, T] $ and all $ M \in \mathbb{N} $. Furthermore, observe that
\begin{align}
\begin{split}
 &\Bigg\|
    \int_{ 0 }^{ t }
    e^{A(t-s)}
    \Bigg[
      F''( \XY{\fl{s}} ) \!
      \left(
        \int_{\fl{s}}^s
        B( \XY{\fl{u}} ) \, dW_u,
        \int_{\fl{s}}^s
        B( \XY{\fl{u}} ) \, dW_u
      \right)
\\&\qquad-
      \sum_{ j \in \mathcal{J} }
      \int_{\fl{s}}^s
      F''( \XY{\fl{s}} )
      \left(
        B( \XY{\fl{u}} ) g_j, B( \XY{\fl{u}} ) g_j
      \right) du
    \Bigg] \, ds
  \Bigg\|_{ L^2( \Omega; H ) }^2
\\&=
  \sum_{ l=0 }^{ M\fl{t}/T }
  \Bigg\|
    \int_{ \frac{lT}{M} }^{ \min\left(\frac{(l+1)T}{M},t\right) }
    e^{A(t-s)}
    \Bigg[
      F''( \XY{\fl{s}} ) \!
      \left(
        \int_{\fl{s}}^s
        B( \XY{\fl{u}} ) \, dW_u,
        \int_{\fl{s}}^s
        B( \XY{\fl{u}} ) \, dW_u
      \right)
\\&\qquad-
      \sum_{ j \in \mathcal{J} }
      \int_{\fl{s}}^s
      F''( \XY{\fl{s}} )
      \left(
        B( \XY{\fl{u}} ) g_j, B( \XY{\fl{u}} ) g_j
      \right) du
    \Bigg] \, ds
  \Bigg\|_{ L^2( \Omega; H ) }^2
\\&\leq
  \frac{ \C T }{ M }
  \sum_{ l=0 }^{ M\fl{t}/T }
  \int_{ \frac{lT}{M} }^{ \min\left(\frac{(l+1)T}{M},t\right) }
  \Bigg[
    \Bigg\|
      F''( \XY{\frac{lT}{M}} ) \!
      \left(
        \int_{\frac{lT}{M}}^s
        B( \XY{\fl{u}} ) \, dW_u,
        \int_{\frac{lT}{M}}^s
        B( \XY{\fl{u}} \, dW_u
      \right)
    \Bigg\|_{ L^2( \Omega; H ) }^2 
\\&\qquad+
    \Bigg\|
      \sum_{ j \in \mathcal{J} }
      \int_{\frac{lT}{M}}^s
      F''( \XY{\frac{lT}{M}} )
      \left(
        B( \XY{\fl{u}} ) g_j, B( \XY{\fl{u}} ) g_j
      \right) du
    \Bigg\|_{ L^2( \Omega; H ) }^2 
  \Bigg] \, ds
\end{split}
\end{align}
and Lemma~7.7 in Da Prato and Zabczyk \cite{dz92} hence implies
\begin{align}
\notag
 &\Bigg\|
    \int_{ 0 }^{ t }
    e^{A(t-s)}
    \bigg[
      F''( \XY{\fl{s}} ) \!
      \left(
        \int_{\fl{s}}^s
        B( \XY{\fl{u}} ) \, dW_u,
        \int_{\fl{s}}^s
        B( \XY{\fl{u}} ) \, dW_u
      \right)
\\&\qquad-
\notag
      \sum_{ j \in \mathcal{J} }
      \int_{\fl{s}}^s
      F''( \XY{\fl{s}} )
      \left(
        B( \XY{\fl{u}} ) g_j, B( \XY{\fl{u}} ) g_j
      \right) du
    \Bigg] \, ds
  \Bigg\|_{ L^2( \Omega; H ) }^2
\\&\leq
\notag
  \frac{ \C }{ M }
  \sum_{ l=0 }^{ M\fl{t}/T }
  \int_{ \frac{lT}{M} }^{ \min\left(\frac{(l+1)T}{M},t\right)  }
  \Bigg[
    \bigg\|
      \int_{\frac{lT}{M}}^s \!
      B( \XY{\fl{u}} ) \, dW_u
    \bigg\|_{ L^4( \Omega; H ) }^4
    +
    \frac{T}{M}
    \int_{\frac{lT}{M}}^s
    \bigg\|
      B( \XY{\fl{u}} )
    \bigg\|_{ L^4( \Omega; HS( U_0, H) ) }^4 du
  \Bigg] \, ds
\\&\leq
\label{eq:tay_martingal}
  \frac{ \C }{ M }
  \sum_{ l=0 }^{ M\fl{t}/T }
  \int_{ \frac{lT}{M} }^{ \min\left(\frac{(l+1)T}{M},t\right) }
  M^{-2} \, ds
  \leq
  \C M^{-3}
\end{align}
for all $ t \in [0, T] $ and all $ M \in \mathbb{N} $. Finally, combining \eqref{eq:tay_F}--\eqref{eq:tay_martingal} yields
\begin{equation}
\label{eq:tay_F_end}
  \left\|
    \int_{ 0 }^{ t }
    e^{A(t-s)}
    \left(
      F( \XY{s} )
      -
      \TF{s}
    \right) ds 
  \right\|_{ L^2( \Omega; H ) } 
  \leq
  \C
  M^{ -\gamma }
\end{equation}
for all $ t \in [0, T] $ and all $ M \in \mathbb{N} $.

\paragraph{Estimation of \eqref{eq:est_taylor_B}}
\label{sec:est_taylor_B}
Similar as in the previous subsection a Taylor expansion of $ B \colon H \rightarrow HS( U_0, H ) $ and the estimate
\begin{multline}
\label{eq:b''_help2}
  \left\|
    \int_0^1
    \left(
      B''\!\left( \XY{\fl{t}} + r(\XY{t} - \XY{\fl{t}}) \right)
      -
      B''( \XY{\fl{t}} )
    \right)
    \left(
      \XY{t} - \XY{\fl{t}}, \XY{t} - \XY{\fl{t}}
    \right)
    \left( 1 - r \right) dr
  \right\|_{ L^2( \Omega; HS( U_0, H ) ) }
\\\leq
  \C M^{ -\frac{3}{2} }
\end{multline}
for all $ t \in [0, T] $ and all $ M \in \mathbb{N} $ give
\begin{align}
\notag
 &\left\|
    B( \XY{t} )
    -
    \TB{t}
  \right\|_{ L^2( \Omega; HS( U_0, H ) ) }
\\&\leq
\notag
  \left\|
    B'( \XY{\fl{t}} )
    \left( 
      \XY{t} - \XY{\fl{t}}
      -
      \int_{\fl{t}}^t
      \left(
        A \XY{\fl{s}}
        +
        F( \XY{\fl{s}} )
      \right) ds
      -
      \int_{\fl{t}}^t
      B( \XY{s} ) \, dW_s
    \right)
  \right\|_{ L^2( \Omega; HS( U_0, H ) ) }
\\&+
\notag
  \left\|
    B'( \XY{\fl{t}} )
    \int_{\fl{t}}^t
    \left(
      B( \XY{s} )
      -
      B( \XY{\fl{s}} )
      -
      B'( \XY{\fl{s}} )
      \int_{\fl{s}}^s
      B( \XY{\fl{u}} ) \, dW_u
    \right) dW_s
  \right\|_{ L^2( \Omega; HS( U_0, H ) ) }
\\&+
\notag
  \left\|
    B''( \XY{\fl{t}} )
    \left( 
      \XY{t} - \XY{\fl{t}}
      -
      \int_{ \fl{t} }^t
      B( \XY{\fl{s}} ) \,  dW_s,
      \XY{t} - \XY{\fl{t}}
      +
      \int_{ \fl{t} }^t
      B( \XY{\fl{s}} ) \,  dW_s
    \right)
  \right\|_{ L^2( \Omega; HS( U_0, H ) ) }
\\&+
\label{eq:tay_B}
  \C M^{ -\frac{3}{2} }
\end{align}
for all $ t \in [0, T] $ and all $ M \in \mathbb{N} $. Moreover, \eqref{eq:det_solution} implies
\begin{align}
\label{eq:tay_B'}
\begin{split}
 &\bigg\|
    B'( \XY{\fl{t}} )
    \left( 
      \XY{t} - \XY{\fl{t}}
      \! - \!
      \int_{\fl{t}}^t \!\!
      \left(
        A \XY{\fl{s}}
        +
        F( \XY{\fl{s}} )
      \right) ds
      \! - \!
      \int_{\fl{t}}^t \!\!
      B( \XY{s} ) \, dW_s
    \right)
  \bigg\|_{ L^2( \Omega; HS( U_0, H ) ) }
\\&\leq
  \C
  \left\|
    \XY{t} - \XY{\fl{t}}
    -
    \int_{\fl{t}}^t
    \left(
      A \XY{\fl{s}}
      +
      F( \XY{\fl{s}} )
    \right) ds
    -
    \int_{\fl{t}}^t
    B( \XY{s} ) \, dW_s
  \right\|_{ L^2( \Omega; H ) }
\\&\leq
  \C M^{ -\gamma }
\end{split}
\end{align}
and a further Taylor expansion of $ B \colon H \rightarrow HS( U_0, H ) $, \eqref{eq:minus_noise} and \eqref{eq:hoelder} show
\begin{align}
\notag
 &\left\|
    B'( \XY{\fl{t}} )
    \int_{\fl{t}}^t
    \left(
      B( \XY{s} )
      -
      B( \XY{\fl{s}} )
      -
      B'( \XY{\fl{s}} )
      \int_{\fl{s}}^s
      B( \XY{\fl{u}} ) \, dW_u
    \right) dW_s
  \right\|_{ L^2( \Omega; HS( U_0, H ) ) }^2
\\&\leq
\label{eq:tay_B'B}
  \C
  \int_{\fl{t}}^t
  \left\|
    B'( \XY{\fl{s}} )
    \left( 
      \XY{s}
      -
      \XY{\fl{s}}
      -
      \int_{\fl{s}}^s
      B( \XY{\fl{u}} ) \, dW_u
    \right)
  \right\|_{ L^2( \Omega; HS( U_0, H ) ) }^2 ds
\\&+
\notag
  \C
  \int_{\fl{t}}^t
  \left\|
    \int_0^1
    B''\!\left( \XY{\fl{s}} + r( \XY{s} - \XY{\fl{s}} ) \right)
    \left( 
      \XY{s}
      -
      \XY{\fl{s}},
      \XY{s}
      -
      \XY{\fl{s}}
    \right)
    ( 1 - r ) \, dr
  \right\|_{ L^2( \Omega; HS( U_0, H ) ) }^2 ds
\\&\leq
\notag
  \C
  \int_{\fl{t}}^t
  \left\|
    \XY{s}
    -
    \XY{\fl{s}}
    -
    \int_{\fl{s}}^s
    B( \XY{\fl{u}} ) \, dW_u
  \right\|_{ L^2( \Omega; H ) }^2 ds
  +
  \C
  \int_{\fl{t}}^t
  \left\|
    \XY{s}
    -
    \XY{\fl{s}}
  \right\|_{ L^4( \Omega; H ) }^4 ds
\\&\leq
\notag
  \C M^{ -3 }
\end{align}
for all $ t \in [0, T] $ and all $ M \in \mathbb{N} $. Furthermore, Lemma~7.7 in Da Prato and Zabczyk \cite{dz92}, \eqref{eq:minus_noise} and \eqref{eq:hoelder} imply
\begin{align}
\label{eq:tay_B''}
\begin{split}
 &\left\|
    B''( \XY{\fl{t}} )
    \left( 
      \XY{t} - \XY{\fl{t}}
      -
      \int_{ \fl{t} }^t
      B( \XY{\fl{s}} ) \,  dW_s,
      \XY{t} - \XY{\fl{t}}
      +
      \int_{ \fl{t} }^t
      B( \XY{\fl{s}} ) \,  dW_s
    \right)
  \right\|_{ L^2( \Omega; HS( U_0, H ) ) }
\\&\leq
  \left\|
    \XY{t}
    -
    \XY{\fl{t}}
    -
    \int_{ \fl{t} }^t
    B( \XY{\fl{s}} ) \,  dW_s
  \right\|_{ L^4( \Omega; H ) }
  \left\|
    \XY{t}
    -
    \XY{\fl{t}}
    +
    \int_{ \fl{t} }^t
    B( \XY{\fl{s}} ) \,  dW_s
  \right\|_{ L^4( \Omega; H ) }
\\&\leq
  \C M^{ -1 }
  \left[
    \left\|
      \XY{t}
      -
      \XY{\fl{t}}
    \right\|_{ L^4( \Omega; H ) }
    +    
    \left\|
      \int_{ \fl{t} }^t
      B( \XY{\fl{s}} ) \,  dW_s
    \right\|_{ L^4( \Omega; H ) }
  \right]
  \leq
  \C M^{ -\frac{3}{2} }
\end{split}
\end{align}
for all $ t \in [0, T] $ and all $ M \in \mathbb{N} $ due to \eqref{eq:minus_noise}. Combining \eqref{eq:tay_B}--\eqref{eq:tay_B''} then finally yields
\begin{equation}
\label{eq:tay_B_end}
  \left\|
    B( \XY{t} )
    -
    \TB{t}
  \right\|_{ L^2( \Omega; HS( U_0, H ) ) }
  \leq 
  \C M^{ -\gamma }
\end{equation}
for all $ t \in [0, T] $ and all $ M \in \mathbb{N} $.

\subsubsection{Estimates for $ \big\| \TZ{m} - \Y{m} \big\|_{ L^2( \Omega; H ) }$ for $ m \in \{ 0,1, \ldots, M \} $ and $ M \in \mathbb{N} $}
The following well known lemma is used in this subsection. 
\begin{lemma}
\label{lem:lemma_eA2}
Let the setting in Section~\ref{sec:sec2} be fulfilled. Then
\begin{equation}
\label{eq:lemm_eA2_a}
  \left\|
    \left(-tA\right)^{-\kappa}
    \left(
      e^{As}
      -
      e^{A\frac{t}{2}}
    \right)
  \right\|_{L(H)}
  \leq
  1
\end{equation}
for all $ s \in [0,t] $, all $ t \in (0, \infty) $ and all $ \kappa \in [0,1] $ and
\begin{equation}
\label{eq:lemm_eA2_b}
 \left\|
    \int_0^t
    \left(-tA\right)^{-\kappa}
    \left(
      e^{Au}
      -
      e^{A\frac{t}{2}}
    \right) du
  \right\|_{L(H)}
  \leq
  t,
\end{equation}
\begin{equation}
\label{eq:lemm_eA2_c}
  \left\|
    \left( -tA \right)^{ -\kappa }
    \left(
      e^{As}
      -
      e^{ A\frac{t}{2} }
      \left(
        I + (s-\tfrac{t}{2})A
      \right) 
    \right)
  \right\|_{ L( H ) }
  \leq
  2
\end{equation}
for all $ s \in [0,t] $, all $ t \in (0, \infty) $ and all $ \kappa \in [0,2] $.
\end{lemma}

\paragraph{Estimation of \eqref{eq:discrete_F}}
\label{sec:discrete_F}
First of all, note that
\begin{align}
\label{eq:dF_start}
\begin{split}
 &\left\|
    \int_{ 0 }^{ \frac{mT}{M} }
    \left( 
      e^{A(\frac{mT}{M}-s)}
      -
      e^{ A(\frac{mT}{M}-\m{s}) }
    \right) 
    \TF{s} \, ds
  \right\|_{ L^2( \Omega; H ) }
\\&\leq
  \left\|
    \int_{ \frac{(m-1)T}{M} }^{ \frac{mT}{M} }
    \left( 
      e^{A(\frac{mT}{M}-s)}
      -
      e^{ A\frac{T}{2M} }
    \right) 
    \TF{s} \, ds
  \right\|_{ L^2( \Omega; H ) }
\\&+
  \sum_{ l=0 }^{ m-2 }
  \left\|
    \left( -A \right)
    e^{A\frac{(m-l-1)T}{M}}
  \right\|_{ L(H) }
  \left\|
    \int_{ \frac{lT}{M} }^{ \frac{(l+1)T}{M} }
    \left( -A \right)^{-1}
    \left( 
      e^{A(\frac{(l+1)T}{M}-s)}
      -
      e^{ A\frac{T}{2M} }
    \right)
    \TF{s} \, ds
  \right\|_{ L^2( \Omega; H ) }
\end{split}
\end{align}
for all $ m \in \{ 1, 2,\ldots, M \} $ and $ M \in \mathbb{N} $. Moreover, observe that
\begin{align}
\notag
  &\left\|
    \int_{ \frac{kT}{M} }^{ \frac{(k+1)T}{M} }
    \left( -A \right)^{ -\kappa }
    \left( 
      e^{A(\frac{(k+1)T}{M}-s)}
      -
      e^{ A\frac{T}{2M} }
    \right) 
    \TF{s} \, ds
  \right\|_{ L^2( \Omega; H ) }
\\&\leq
  \left\|
    \int_{ \frac{kT}{M} }^{ \frac{(k+1)T}{M} }
    \left( -A \right)^{ -(\kappa+\alpha) }
    \left(
      e^{A(\frac{(k+1)T}{M}-s)}
      -
      e^{ A\frac{T}{2M} }
    \right) ds
  \right\|_{ L(H) }
  \left\|
    F( \Y{k} )
  \right\|_{ L^2( \Omega; H_{ \alpha } ) }
\\&+
\notag
  \int_{ \frac{kT}{M} }^{ \frac{(k+1)T}{M} }
  \left\|
    \left( -A \right)^{ -\kappa }
    \left(
      e^{A(\frac{(k+1)T}{M}-s)}
      -
      e^{ A\frac{T}{2M} }
    \right)
  \right\|_{ L(H) }
  \left\|
    F'(\Y{k})
    \int_{\frac{kT}{M}}^s
    \left(
      A \Y{k}
      +
      F(\Y{k})
    \right) du
  \right\|_{ L^2( \Omega; H ) } ds
\\&+
\notag
  \int_{ \frac{kT}{M} }^{ \frac{(k+1)T}{M} }
  \left\|
    \left( -A \right)^{ -\kappa }
    \left(
      e^{A(\frac{(k+1)T}{M}-s)}
      -
      e^{ A\frac{T}{2M} }
    \right)
  \right\|_{ L(H) }
  \left\|
    F'(\Y{k})
    \int_{\frac{kT}{M}}^s
    B(\Y{k}) \, dW_u
  \right\|_{ L^2( \Omega; H ) } ds
\\&+
\notag
  \int_{ \frac{kT}{M} }^{ \frac{(k+1)T}{M} }
  \int_{ \frac{kT}{M} }^s
  \left\|
    \left( -A \right)^{ -\kappa } 
    \left( 
      e^{A(\frac{(k+1)T}{M}-s)}
      -
      e^{ A\frac{T}{2M} }
    \right)
    \sum_{ j \in \mathcal{J} }
    F''( \Y{k} )
    \left(
      B( \Y{k} ) g_j, B( \Y{k} ) g_j
    \right) 
  \right\|_{ L^2( \Omega; H ) } du \, ds
\end{align}
and Lemma~\ref{lem:lemma_eA2} hence gives 
\begin{align}
\begin{split}
  &\left\|
    \int_{ \frac{kT}{M} }^{ \frac{(k+1)T}{M} }
    \left( -A \right)^{ -\kappa }
    \left( 
      e^{A(\frac{(k+1)T}{M}-s)}
      -
      e^{ A\frac{T}{2M} }
    \right) 
    \TF{s} \, ds
  \right\|_{ L^2( \Omega; H ) }
\\&\leq
\label{eq:help1}
  \C M^{ -(1 + \min(\kappa+\alpha,2)) }
  +
  \C M^{ -\kappa }
  \int_{ \frac{kT}{M} }^{ \frac{(k+1)T}{M} }
  \int_{ \frac{kT}{M} }^s
  \left\|
    A \Y{k}
    +
    F(\Y{k})
  \right\|_{ L^2( \Omega; H ) } du \, ds
\\&+
  \C M^{ -\kappa }
  \int_{ \frac{kT}{M} }^{ \frac{(k+1)T}{M} }
  \left[ 
    \int_{ \frac{kT}{M} }^s
    \left\|
      B(\Y{k})
    \right\|_{ L^2( \Omega; HS( U_0, H ) ) }^2 du 
  \right]^{ \frac{1}{2} } ds
\\&+
  \C M^{-(\kappa+1)}
  \int_{ \frac{kT}{M} }^{ \frac{(k+1)T}{M} }
  \left\|
    \sum_{ j \in \mathcal{J} }
    F''(\Y{k})
    \left(
      B(\Y{k}) g_j, B(\Y{k}) g_j
    \right) 
  \right\|_{ L^2( \Omega; H ) } ds
\\&\leq
  \C M^{ -(\kappa + 1 +\min(\alpha, \frac{1}{2})) }
\end{split}
\end{align}
for all $ k \in \{ 0,1,\ldots,M-1\} $, all $ M \in \mathbb{N} $ and all $ \kappa \in [0,1] $. Combining \eqref{eq:dF_start} and \eqref{eq:help1} implies
\begin{align}
\notag
 &\left\|
    \int_{ 0 }^{ \frac{mT}{M} }
    \left( 
      e^{A(\frac{mT}{M}-s)}
      -
      e^{ A(\frac{mT}{M}-\m{s}) }
    \right) 
    \TF{s} \, ds
  \right\|_{ L^2( \Omega; H ) }
\\&\leq
\label{eq:discrete_F_end}
  \C M^{ -(1 + \min(\alpha, \frac{1}{2})) }
  +
  \C M^{ -(2 + \min(\alpha, \frac{1}{2})) }
  \sum_{ l=1 }^{ m-1 }
  \left\|
    \left( -A \right)
    e^{A\frac{lT}{M}}
  \right\|_{ L( H ) }
\\&\leq
\notag
  \C M^{ -(1 + \min(\alpha, \frac{1}{2})) }
  +
  \C M^{ -(1 + \min(\alpha, \frac{1}{2})) }
  \left[ 
    \sum_{ l=1 }^{ m-1 }
    \frac{1}{l}
  \right]
  \leq
  \C M^{ -\min(\alpha+1,\frac{3}{2}) }
  \left(
    1
    +
    \log(M)
  \right)
\end{align}
for all $ m \in \{ 0,1,\ldots,M\} $ and all $ M \in \mathbb{N} $.

\paragraph{Estimation of \eqref{eq:discrete_B}}
\label{sec:discrete_B}
First of all, observe that
\begin{align}
\begin{split}
 &\left\|
    \int_{ 0 }^{ \frac{mT}{M} }
    \left( 
      e^{A(\frac{mT}{M}-s)} \,
      \TB{s}
      -
      e^{A(\frac{mT}{M}-\m{s})}
      \left(
        \TB{s}
        +
        \left( \m{s} \!-\! s \right)\! A
        B( \XY{\fl{s}} )
      \right)
    \right) dW_s
  \right\|_{ L^2( \Omega; H ) }^2
\\&\leq
  \int_{ \frac{(m-1)T}{M} }^{ \frac{mT}{M} }  
  \left\| 
      e^{A(\frac{mT}{M}-s)} \,
      \TB{s}
      -
      e^{ A\frac{T}{2M} }
      \left(
        \TB{s}
        +
        \left( \m{s} \!-\! s \right)\! A
        B( \Y{(m-1)} )
      \right)
  \right\|_{ L^2( \Omega; HS( U_0, H ) ) }^2 ds
\\&+
\label{eq:discrete_B_s}
  \sum_{ l=0 }^{ m-2 }
  \int_{ \frac{lT}{M} }^{ \frac{(l+1)T}{M} }
  \left\|
    \left(-A\right)^{ \frac{1}{2} }
    e^{ A\frac{(m-l-1)T}{M} }
  \right\|^2
\\&\qquad\qquad\cdot
  \left\|
    \left(-A\right)^{ -\frac{1}{2} }
    \left( 
      e^{A(\frac{(l+1)T}{M}-s)} \,
      \TB{s}
      -
      e^{ A\frac{T}{2M} }
      \left(
        \TB{s}
        +
        \left( \m{s} \!-\! s \right)\! A
        B( \Y{l} )
      \right)
    \right)
  \right\|_{ L^2( \Omega; HS( U_0, H ) ) }^2 ds
\end{split}
\end{align}
for all $ m \in \{ 1,2,\ldots,M\} $ and all $ M \in \mathbb{N} $. In addition, note that
\begin{align}
\begin{split}
 &\left\|
    \left( -A \right)^{ -\kappa }
    \left(
      e^{A(\frac{(k+1)T}{M}-s)} \,
      \TB{s}
      -
      e^{ A\frac{T}{2M} }
      \left(
        \TB{s}
        +
        \left( \m{s} \!-\! s \right)\! A
        B( \Y{k} )
      \right)
    \right)
  \right\|_{ L^2( \Omega; HS( U_0, H ) ) }
\\&\leq
  \left\|
    \left( -A \right)^{ -(\kappa+\beta) }
    \left( 
      e^{A(\frac{(k+1)T}{M}-s)}
      -
      e^{ A\frac{T}{2M} }
      \left(
        I + \left( \m{s} \!-\! s \right)\! A
      \right) 
    \right)
  \right\|_{ L( H ) }
  \left\|
    \left( -A \right)^{ \beta }
    B( \Y{k} ) 
  \right\|_{ L^2( \Omega; HS( U_0, H ) ) }
\\&+
  \left\|
    \left( -A \right)^{ -\kappa }
    \left( 
      e^{A(\frac{(k+1)T}{M}-s)}
      -
      e^{ A\frac{T}{2M} }
    \right)
  \right\|_{ L( H ) }
  \left\|
    B'( \Y{k} )
    \int_{ \frac{kT}{M} }^s
    \left( 
      A \Y{k}
      +
      F( \Y{k} ) 
    \right) du
  \right\|_{ L^2( \Omega; HS( U_0, H ) ) }
\\&+
  \left\|
    \left( -A \right)^{ -(\kappa+\delta) }
    \left( 
      e^{A(\frac{(k+1)T}{M}-s)}
      -
      e^{ A\frac{T}{2M} }
    \right)
  \right\|_{ L( H ) }
  \left\|
    \left( -A \right)^{ \delta }
    B'( \Y{k} )
    \int_{ \frac{kT}{M} }^s
    B( \Y{k} ) \, dW_u
  \right\|_{ L^2( \Omega; HS( U_0, H ) ) }
\\&+
  \C
  \left\|
    \left( -A \right)^{ -\kappa }
    \left( 
      e^{A(\frac{(k+1)T}{M}-s)}
      -
      e^{ A\frac{T}{2M} }
    \right)
  \right\|_{ L( H ) }
  \left\|
    \int_{ \frac{kT}{M} }^s
    B'( \Y{k} )
    \left(
      \int_{ \frac{kT}{M} }^u
      B( \Y{k} ) \, dW_v
    \right) dW_u
  \right\|_{ L^2( \Omega; H ) }
\\&+
  \C
  \left\|
    \left( -A \right)^{ -\kappa }
    \left( 
      e^{A(\frac{(k+1)T}{M}-s)}
      -
      e^{ A\frac{T}{2M} }
    \right)
  \right\|_{ L( H ) }
  \left\|
    \int_{ \frac{kT}{M} }^s
    B( \Y{k} ) \, dW_u
  \right\|_{ L^4( \Omega; H ) }^2
\end{split}
\end{align}
and Lemma~7.7 in Da Prato and Zabczyk \cite{dz92} and Lemma~\ref{lem:lemma_eA2} hence imply
\begin{align}
\notag
 &\left\|
    \left( -A \right)^{ -\kappa }
    \left(
      e^{A(\frac{(k+1)T}{M}-s)} \,
      \TB{s}
      -
      e^{ A\frac{T}{2M} }
      \left(
        \TB{s}
        +
        \left( \m{s} \!-\! s \right)\! A
        B( \Y{k} )
      \right)
    \right)
  \right\|_{ L^2( \Omega; HS( U_0, H ) ) }
\\&\leq
\notag
  \C M^{ -(\kappa+\beta) }
  +
  \C M^{ -\kappa }
  \int_{ \frac{kT}{M} }^s
  \left(
    \left\|
      A \Y{k}
    \right\|_{ L^2( \Omega; H ) } 
    +
    \left\|
      F( \Y{k} )
    \right\|_{ L^2( \Omega; H ) } 
  \right) du
\\&+
\label{eq:d_B_2}
  \C M^{ -(\kappa+\delta) }
  \left(
    \int_{ \frac{kT}{M} }^s
    \left\|
      \left( -A \right)^{ \delta }
      B( \Y{k} )
    \right\|_{ L^2( \Omega; HS( U_0, H ) ) }^2 du
  \right)^{ \frac{1}{2} }
\\&+
\notag
  \C M^{ -\kappa }
  \left( 
    \int_{ \frac{kT}{M} }^s
    \int_{ \frac{kT}{M} }^u
    \left\|
      B( \Y{k} )
    \right\|_{ L^2( \Omega; HS( U_0, H ) ) }^2 dv \, du
  \right)^{ \frac{1}{2} }
  +
  \C M^{ -\kappa }
  \int_{ \frac{kT}{M} }^s
  \left\|
    B( \Y{k} )
  \right\|_{ L^4( \Omega; HS( U_0, H ) ) }^2 du
\\&\leq
\notag
  \C M^{ -(\kappa+\beta) }
  +
  \C M^{ -(\kappa+1) }
  +
  \C M^{ -(\kappa+\delta+\frac{1}{2}) }
  \leq
  \C M^{ -(\kappa+\min(\beta, \delta+\frac{1}{2}, 1)) }
\end{align}
for all $ s \in [\frac{kT}{M}, \frac{(k+1)T}{M}] $, all $ k \in \{ 0,1,\ldots,M-1\} $, all $ M \in \mathbb{N} $ and all $ \kappa \in [0, \frac{1}{2}] $. Combining \eqref{eq:discrete_B_s} and \eqref{eq:d_B_2} then shows
\begin{align}
\label{eq:discrete_B_end}
\begin{split}
 &\left\|
    \int_{ 0 }^{ \frac{mT}{M} }
    \left( 
      e^{A(\frac{mT}{M}-s)} \,
      \TB{s}
      -
      e^{A(\frac{mT}{M}-\m{s})}
      \left(
        \TB{s}
        +
        \left( \m{s} \!-\! s \right)\! A
        B( \XY{\fl{s}} )
      \right)
    \right) dW_s
  \right\|_{ L^2( \Omega; H ) }^2
\\&\leq
  \C M^{ -2\min(\beta+\frac{1}{2},\delta+1,\frac{3}{2}) }
  +
  \C M^{ -2\min(\beta+\frac{1}{2},\delta+1,\frac{3}{2}) }
  \left[
    \sum_{ l=1 }^{ m-1 }
    \frac{1}{l}
  \right]
\\&\leq
  \C M^{ -2\min(\beta+\frac{1}{2},\delta+1,\frac{3}{2}) }
  \left( 1 + \log(M) \right)
\end{split}
\end{align}
for all $ m \in \{ 0,1,\ldots,M\} $ and all $ M \in \mathbb{N} $.

\subsection{Proof of Lemma~\ref{lem:noise2}}
\label{sec:proof_of_lemnoise2}

Throughout this proof let $ \I{i}{t_0,t}, \I{i,j}{t_0,t}, \I{i,j,k}{t_0,t} \colon \Omega \rightarrow \mathbb{R} $, $ i,j,k \in \mathcal{J} $, $ t_0, t \in [0,T], t_0 \leq t $, be random variables satisfying
\begin{align}
\begin{split}
  \I{i}{t_0,t}
&=
  \int_{t_0}^t
  \left< g_i, dW_s \right>_{ U_0 }
\\
  \I{i,j}{t_0,t}
&=
  \int_{t_0}^t
  \int_{t_0}^s
  \left< g_i, dW_u \right>_{ U_0 }
  \left< g_j, dW_s \right>_{ U_0 }
\\
  \I{i,j,k}{t_0,t}
&=
  \int_{t_0}^t
  \int_{t_0}^s
  \int_{t_0}^u
  \left< g_i, dW_v \right>_{ U_0 }
  \left< g_j, dW_u \right>_{ U_0 }
  \left< g_k, dW_s \right>_{ U_0 }
\end{split}
\end{align}
$ \mathbb{P} $-a.s. for all $ i,j,k \in \mathcal{J} $ and all $ t_0, t \in [0,T] $ with $ t_0 \leq t $. In addition, we use the following well known identities for stochastic integrals (see, e.g., (10.3.15) and (10.4.14) in \cite{kp92})
\begin{equation}
\label{eq:ito_a}
  \I{i}{t_0,t} \,
  \I{j}{t_0,t}
  =
  \I{i,j}{t_0,t}
  +
  \I{j,i}{t_0,t},
\qquad
  \I{i,i,i}{t_0,t}
  =
  \frac{1}{6}
  \left(
    \I{i}{t_0,t}^2
    -
    3 \left( t - t_0 \right)
  \right)
  \I{i}{t_0,t}
\end{equation}
$ \mathbb{P} $-a.s. for all $ i, j \in \mathcal{J} $ and all $ t_0, t \in [0,T] $ with $ t_0 \leq t $ and
\begin{align}
\label{eq:ito_b}
\begin{split}
 &\I{i,j,k}{t_0,t}
  +
  \I{j,i,k}{t_0,t}
  +
  \I{j,k,i}{t_0,t}
  +
  \I{i,k,j}{t_0,t}
  +
  \I{k,j,i}{t_0,t}
  +
  \I{k,i,j}{t_0,t}
  =
  \I{i}{t_0,t} \,
  \I{j}{t_0,t} \,
  \I{k}{t_0,t}
\\&
  \I{i,j,j}{t_0,t}
  +
  \I{j,i,j}{t_0,t}
  +
  \I{j,j,i}{t_0,t}
  =
  \frac{1}{2}
  \I{i}{t_0,t}
  \left(
    \I{j}{t_0,t}^2
    -
    \left( t - t_0 \right)
  \right)
\end{split}
\end{align}
$ \mathbb{P} $-a.s. for all $ i, j, k \in \mathcal{J}, i \neq j, i \neq k, j \neq k $ and all $ t_0, t \in [0,T] $ with $ t_0 \leq t $. Next note that
\begin{align}
\label{eq:b'b'b_help}
\begin{split}
 &\int_{ t_0 }^{ t }
  B'( Z )
  \left( 
    \int_{ t_0 }^s
    B'( Z ) 
    \left( 
      \int_{ t_0 }^u
      B( Z ) \, dW_v 
    \right) dW_u 
  \right) dW_s
\\&=
  \sum_{ i,j,k \in \mathcal{J} }
  \int_{ t_0 }^{ t }
  B'( Z )
  \left( 
    \int_{ t_0 }^s
    B'( Z ) 
    \left( 
      \int_{ t_0 }^u
      B( Z ) \, g_i \left< g_i, dW_v \right>_{ U_0 }
    \right) g_j \left< g_j, dW_u \right>_{ U_0 }
  \right) g_k \left< g_k, dW_s \right>_{ U_0 }
\\&=
  \sum_{ i,j,k \in \mathcal{J} }
  B'( Z )
  \Big( 
    B'( Z ) 
    \Big( 
      B( Z ) \, g_i 
    \Big) g_j 
  \Big) g_k \,
  \I{i,j,k}{t_0,t}
\end{split}
\end{align}
$ \mathbb{P} $-a.s. for all $ \mathcal{F}_{t_0} / \mathcal{B}(H) $-measurable mappings $ Z \colon \Omega \rightarrow H $ and all $ t_0, t \in [0,T] $ with $ t_0 \leq t $. In addition, \eqref{eq:ito_a} shows
\begin{align}
\label{eq:b''bb_help}
\begin{split}
 &\frac{1}{2}
  \int_{ t_0 }^{ t }
  B''( Z )
  \left(
    \int_{ t_0 }^s
    B( Z ) \, dW_u,
    \int_{ t_0 }^s
    B( Z ) \, dW_u
  \right) dW_s
\\&=
  \frac{1}{2}
  \sum_{ i,j,k \in \mathcal{J} }
  \int_{ t_0 }^{ t }
  B''( Z )
  \left(
    \int_{ t_0 }^s
    B( Z ) \, g_i \left< g_i, dW_u \right>_{ U_0 },
    \int_{ t_0 }^s
    B( Z ) \, g_j \left< g_j, dW_u \right>_{ U_0 }
  \right)
  g_k \, \left< g_k, dW_s \right>_{ U_0 }
\\&=
  \frac{1}{2}
  \sum_{ i,j,k \in \mathcal{J} }
  B''( Z )
  \Big(
    B( Z ) \, g_i,
    B( Z ) \, g_j
  \Big) g_k
  \int_{ t_0 }^{ t }
  \I{i}{t_0,s} \,
  \I{j}{t_0,s} \,
  d\!\left< g_k, W_s \right>_{ U_0 }
\\&=
  \frac{1}{2}
  \sum_{ i,j,k \in \mathcal{J} }
  B''( Z )
  \Big(
    B( Z ) \, g_i,
    B( Z ) \, g_j
  \Big) g_k
  \left( 
    \I{i,j,k}{t_0,t}
    +
    \I{j,i,k}{t_0,t}
  \right)
\\&=
  \sum_{ i,j,k \in \mathcal{J} }
  B''( Z )
  \Big(
    B( Z ) \, g_i,
    B( Z ) \, g_j
  \Big) g_k \,
  \I{i,j,k}{t_0,t}
\end{split}
\end{align}
$ \mathbb{P} $-a.s. for all $ \mathcal{F}_{t_0} / \mathcal{B}(H) $-measurable mappings $ Z \colon \Omega \rightarrow H $ and all $ t_0, t \in [0,T] $ with $ t_0 \leq t $. Combining \eqref{eq:b'b'b_help} and \eqref{eq:b''bb_help} then gives
\begin{align}
\notag
 &\int_{ t_0 }^{ t }
  B'( Z )
  \left( 
    \int_{ t_0 }^s
    B'( Z ) 
    \left( 
      \int_{ t_0 }^u
      B( Z ) \, dW_v 
    \right) dW_u 
  \right) dW_s
  +
  \frac{1}{2}
  \int_{ t_0 }^{ t }
  B''( Z )
  \left(
    \int_{ t_0 }^s
    B( Z ) \, dW_u,
    \int_{ t_0 }^s
    B( Z ) \, dW_u
  \right) dW_s
\\&=
\notag
  \sum_{ i,j,k \in \mathcal{J} }
  \left[ 
  B'( Z )
  \Big( 
    B'( Z ) 
    \Big( 
      B( Z ) \, g_i 
    \Big) g_j 
  \Big)
  +
  B''( Z )
  \Big(
    B( Z ) \, g_i,
    B( Z ) \, g_j
  \Big) 
  \right] g_k \,
  \I{i,j,k}{t_0,t}
\\&=
\notag
  \frac{1}{6}
  \sum_{ \substack{ i,j,k \in \mathcal{J} \\ i \neq j, i \neq k, j \neq k } }
  \left[ 
  B'( Z )
  \Big( 
    B'( Z ) 
    \Big( 
      B( Z ) \, g_i 
    \Big) g_j 
  \Big)
  +
  B''( Z )
  \Big(
    B( Z ) \, g_i,
    B( Z ) \, g_j
  \Big) 
  \right] g_k \,
  \bigg[
    \I{i,j,k}{t_0,t}
    +
    \I{j,i,k}{t_0,t}
\\&\quad+
    \I{j,k,i}{t_0,t}
    +
    \I{i,k,j}{t_0,t}
    +
    \I{k,j,i}{t_0,t}
    +
    \I{k,i,j}{t_0,t}
  \bigg]
\\&+
\notag
  \sum_{ \substack{ i,j \in \mathcal{J} \\ i \neq j } }
  \left[ 
  B'( Z )
  \Big( 
    B'( Z ) 
    \Big( 
      B( Z ) \, g_i 
    \Big) g_i
  \Big)
  \!+\!
  B''( Z )
  \Big(
    B( Z ) \, g_i,
    B( Z ) \, g_i
  \Big) 
  \right] g_j
  \left[ 
    \I{i,j,i}{t_0,t}
    +
    \I{j,i,i}{t_0,t}
    +
    \I{i,i,j}{t_0,t}
  \right]
\\&+
\notag
  \sum_{ i \in \mathcal{J} }
  \left[ 
  B'( Z )
  \Big( 
    B'( Z ) 
    \Big( 
      B( Z ) \, g_i 
    \Big) g_i 
  \Big)
  +
  B''( Z )
  \Big(
    B( Z ) \, g_i,
    B( Z ) \, g_i
  \Big) 
  \right] g_i \,
  \I{i,i,i}{t_0,t}
\end{align}
and \eqref{eq:ito_a}--\eqref{eq:ito_b} hence show
\begin{align}
\notag
 &\int_{ t_0 }^{ t }
  B'( Z )
  \left( 
    \int_{ t_0 }^s
    B'( Z ) 
    \left( 
      \int_{ t_0 }^u
      B( Z ) \, dW_v 
    \right) dW_u 
  \right) dW_s
  +
  \frac{1}{2}
  \int_{ t_0 }^{ t }
  B''( Z )
  \left(
    \int_{ t_0 }^s
    B( Z ) \, dW_u,
    \int_{ t_0 }^s
    B( Z ) \, dW_u
  \right) dW_s
\\&=
\notag
  \frac{1}{6}
  \sum_{ \substack{ i,j,k \in \mathcal{J} \\ i \neq j, i \neq k, j \neq k } }
  \left[ 
  B'( Z )
  \Big( 
    B'( Z ) 
    \Big( 
      B( Z ) \, g_i 
    \Big) g_j 
  \Big)
  +
  B''( Z )
  \Big(
    B( Z ) \, g_i,
    B( Z ) \, g_j
  \Big) 
  \right] g_k \,
  \I{i}{t_0,t} \,
  \I{j}{t_0,t} \,
  \I{k}{t_0,t}
\\&+
  \frac{1}{2}
  \sum_{ \substack{ i,j \in \mathcal{J} \\ i \neq j } }
  \left[ 
  B'( Z )
  \Big( 
    B'( Z ) 
    \Big( 
      B( Z ) \, g_i 
    \Big) g_i
  \Big)
  +
  B''( Z )
  \Big(
    B( Z ) \, g_i,
    B( Z ) \, g_i
  \Big) 
  \right] g_j
  \I{j}{t_0,t}
  \left(
    \I{i}{t_0,t}^2
    -
    \left( t - t_0 \right)
  \right)
\\&+
\notag
  \frac{1}{6}
  \sum_{ i \in \mathcal{J} }
  \left[ 
  B'( Z )
  \Big( 
    B'( Z ) 
    \Big( 
      B( Z ) \, g_i 
    \Big) g_i 
  \Big)
  +
  B''( Z )
  \Big(
    B( Z ) \, g_i,
    B( Z ) \, g_i
  \Big) 
  \right] g_i \,
  \left(
    \I{i}{t_0,t}^2
    -
    3 \left( t - t_0 \right)
  \right)
  \I{i}{t_0,t}
\end{align}
$ \mathbb{P} $-a.s. for all $ \mathcal{F}_{t_0} / \mathcal{B}(H) $-measurable mappings $ Z \colon \Omega \rightarrow H $ and all $ t_0, t \in [0,T] $ with $ t_0 \leq t $. This finally yields
\begin{align}
\notag
 &\int_{ t_0 }^{ t }
  B'( Z )
  \left( 
    \int_{ t_0 }^s
    B'( Z ) 
    \left( 
      \int_{ t_0 }^u
      B( Z ) \, dW_v 
    \right) dW_u 
  \right) dW_s
  +
  \frac{1}{2}
  \int_{ t_0 }^{ t }
  B''( Z )
  \left(
    \int_{ t_0 }^s
    B( Z ) \, dW_u,
    \int_{ t_0 }^s
    B( Z ) \, dW_u
  \right) dW_s
\\&=
\notag
  \frac{1}{6}
  \sum_{ i,j,k \in \mathcal{J} }
  \left[ 
  B'( Z )
  \Big( 
    B'( Z ) 
    \Big( 
      B( Z ) \, g_i 
    \Big) g_j 
  \Big)
  +
  B''( Z )
  \Big(
    B( Z ) \, g_i,
    B( Z ) \, g_j
  \Big) 
  \right] g_k \,
  \I{i}{t_0,t} \,
  \I{j}{t_0,t} \,
  \I{k}{t_0,t}
\\&-
\notag
  \frac{\left( t - t_0 \right)}{2}
  \sum_{ i,j \in \mathcal{J} }
  \left[ 
  B'( Z )
  \Big( 
    B'( Z ) 
    \Big( 
      B( Z ) \, g_i 
    \Big) g_i 
  \Big)
  +
  B''( Z )
  \Big(
    B( Z ) \, g_i,
    B( Z ) \, g_i
  \Big) 
  \right] g_j \,
  \I{j}{t_0,t}
\\&=
  \frac{1}{6}
  B'( Z )
  \Big( 
    B'( Z ) 
    \Big( 
      B( Z ) \left(W_t - W_{t_0}\right)
    \Big) \left(W_t - W_{t_0}\right)
  \Big)
  \left(W_t - W_{t_0}\right)
\\&+
\notag
  \frac{1}{6}
  B''( Z )
  \Big(
    B( Z ) \left(W_t - W_{t_0}\right),
    B( Z ) \left(W_t - W_{t_0}\right)
  \Big) 
  \left(W_t - W_{t_0}\right)
\\&-
\notag
  \frac{\left( t - t_0 \right)}{2}
  \sum_{ i \in \mathcal{J} }
  \left[ 
  B'( Z )
  \Big( 
    B'( Z ) 
    \Big( 
      B( Z ) \, g_i 
    \Big) g_i 
  \Big)
  \left(W_t - W_{t_0}\right)
  +
  B''( Z )
  \Big(
    B( Z ) \, g_i,
    B( Z ) \, g_i
  \Big) 
  \left(W_t - W_{t_0}\right)
  \right]
\end{align}
$ \mathbb{P} $-a.s. for all $ \mathcal{F}_{t_0} / \mathcal{B}(H) $-measurable mappings $ Z \colon \Omega \rightarrow H $ and all $ t_0, t \in [0,T] $ with $ t_0 \leq t $.

\subsection{Proof of Lemma~\ref{lem:int}}
\label{sec:proof_of_lemint}

First of all, note that
\begin{align}
\label{eq:int_h1}
\begin{split}
 &\mathbb{E}
   \!\left[
    \big< v, \int_{t_0}^{t} ( W_{s_1} - W_{t_0} ) \, ds_1 \big>_U
    \big< w, \int_{t_0}^{t} ( W_{s_2} - W_{t_0} ) \, ds_2 \big>_U
  \right]
\\&=
  \int_{t_0}^{t}
  \int_{t_0}^{t}
  \mathbb{E}
    \!\left[
      \big< v, ( W_{s_1} - W_{t_0} ) \big>_U
      \big< w, ( W_{s_2} - W_{t_0} ) \big>_U
  \right] ds_1 \, ds_2
\\&=
  \int_{t_0}^{t}
  \int_{t_0}^{t}
  \big< v,\big(\min(s_1,s_2) - t_0 \big) Q w \big>_U \, ds_1 \, ds_2
\\&=
  \int_{t_0}^{t}
  \int_{t_0}^{s_2}
  \big< v,s_1 Q w \big>_U \, ds_1 \, ds_2
  +
  \int_{t_0}^{t}
  \int_{s_2}^{t}
  \big< v,s_2 Q w \big>_U \, ds_1 \, ds_2
  -
  \int_{t_0}^{t}
  \int_{t_0}^{t}
  \big< v,t_0 Q w \big>_U \, ds_1 \, ds_2
\\&=
  \big< v,\tfrac{1}{3}(t - t_0)^3 Q w \big>_U
\end{split}
\end{align}
for all $ v, w \in U $ and all $ t_0, t \in [0,T] $ with $ t_0 \leq t $. Moreover, observe that
\begin{align}
\label{eq:int_h2}
\begin{split}
 &\mathbb{E}
   \!\left[
    \big< v, ( W_{t} - W_{t_0} ) \big>_U
    \big< w, \int_{t_0}^{t} ( W_{s} - W_{t_0} ) \, ds \big>_U
  \right]
\\&=
  \int_{t_0}^{t}
  \mathbb{E}
    \!\left[
      \big< v, ( W_{t} - W_{t_0} ) \big>_U
      \big< w, ( W_{s} - W_{t_0} ) \big>_U
  \right] ds
  =
  \int_{t_0}^{t}
  \big< v,\big(s - t_0 \big) Q w \big>_U \, ds
\\&=
  \big< v,\tfrac{1}{2}(t - t_0)^2 Q w \big>_U
\end{split}
\end{align}
for all $ v, w \in U $ and all $ t_0, t \in [0,T] $ with $ t_0 \leq t $. Combining \eqref{eq:int_h1} and \eqref{eq:int_h2} then yields
\begin{align}
\begin{split}
 &\mathbb{E}
  \!\left[ 
    \left(
      \begin{array}{c}
        \big< v, ( W_{t} - W_{t_0} ) \big>_U \\
        \big< v, \int_{t_0}^{t} ( W_{s} - W_{t_0} ) \, ds \big>_U
      \end{array}
    \right)
    \left(
      \begin{array}{c}
        \big< w, ( W_{t} - W_{t_0} ) \big>_U \\
        \big< w, \int_{t_0}^{t} ( W_{s} - W_{t_0} ) \, ds \big>_U
      \end{array}
    \right)^{\!\!*\,} 
  \right] 
\\&=
  \left(
   \begin{array}{l l}
     \big< v,(t - t_0) Q w \big>_U 
    &\big< v,\tfrac{1}{2}(t - t_0)^2 Q w \big>_U
     \\
     \big< v,\tfrac{1}{2}(t - t_0)^2 Q w\big>_U
    &\big< v,\tfrac{1}{3}(t - t_0)^3 Q w\big>_U
   \end{array}
 \right)
\end{split}
\end{align}
for all $ v, w \in U $ and all $ t_0, t \in [0,T] $ with $ t_0 \leq t $.

\bibliographystyle{acm}
\bibliography{bibfile}

\end{document}